\documentclass[10pt,journal,twocolumn,twoside]{IEEEtran} 

 \pdfoutput=1
\usepackage{amsmath, amssymb}
\usepackage{amsfonts}
\usepackage{graphicx}
\usepackage{enumerate}
\usepackage[font = small]{caption}
\usepackage{ifthen}
\usepackage{multicol}
\usepackage{float}
\usepackage{cite}
\usepackage{fancyhdr}
\usepackage[usenames, dvipsnames]{color}
\usepackage[lofdepth,lotdepth]{subfig}
\graphicspath{{figures/}}
\usepackage{mathtools}

 \usepackage{tikz}
 \usetikzlibrary{decorations.pathmorphing,shadings}
\usepackage{tkz-graph}
\usetikzlibrary{arrows}

 \usepackage{cancel} %To be able to cross out text

\newtheorem{theorem}{Theorem}[section]
\newtheorem{lemma}[theorem]{Lemma}

\newtheorem{corollary}[theorem]{Corollary}
\newtheorem{rem}{Remark}
\newtheorem{example}{Example}
\newtheorem{ass}{Assumption}
\newtheorem{definition}{Definition}

\newboolean{showcomments}
\setboolean{showcomments}{true}

\newcommand{\todo}[1]{  \ifthenelse{\boolean{showcomments}}
{\textcolor{ForestGreen}{TO DO:  #1}}{}}
\newcommand{\bassam}[1]{\ifthenelse{\boolean{showcomments}}
{\textcolor{Orange}{(Bassam says: #1)}}{}}
\newcommand{\henrik}[1]{\ifthenelse{\boolean{showcomments}}
{\textcolor{Blue}{(Henrik says: #1)}}{}}
\newcommand{\partha}[1]{\ifthenelse{\boolean{showcomments}}
{\textcolor{Blue}{(Partha says: #1)}}{}}
\newcommand{\emma}[1]{\ifthenelse{\boolean{showcomments}}
{\textcolor{VioletRed}{(Emma says: #1)}}{}}
\newcommand{\ifneeded}[1]{\ifthenelse{\boolean{showcomments}}
{\textcolor{Gray}{#1}}{}}

\newboolean{showedit}
\setboolean{showedit}{true}
\newcommand{\edit}[1]{\ifthenelse{\boolean{showedit}}
{\textcolor{Blue}{#1}}{}}

\newcommand{\hn}{$\mathcal{H}_2$ }
\newcommand{\hnnorm}{$\mathcal{H}_2$-norm }

\newcommand{\Znd}{\mathbb{Z}_L^d}
\newcommand{\xtr}{\hat{p}(\theta)}
\newcommand{\Vn}{V_N}

\usepackage{array}
\newcolumntype{L}[1]{>{\raggedright\let\newline\\\arraybackslash\hspace{0pt}}m{#1}}
\newcolumntype{C}[1]{>{\centering\let\newline\\\arraybackslash\hspace{0pt}}m{#1}}
\newcolumntype{R}[1]{>{\raggedleft\let\newline\\\arraybackslash\hspace{0pt}}m{#1}}
\usepackage{subfig}

\usepackage{tikz}
 \usetikzlibrary{plotmarks}
 \usepackage{pgfplots}

\usetikzlibrary{circuits.ee.IEC, positioning, external}
\definecolor{gray3}{rgb}{0.75, 0.75, 0.75}
\definecolor{gray2}{rgb}{0.5, 0.5, 0.5}
\definecolor{gray1}{rgb}{0.25, 0.25, 0.25}
\definecolor{gray0}{rgb}{0.15, 0.15, 0.15}

\definecolor{emmagreen1}{rgb}{0, 0.5, 0.1}
\definecolor{emmaorange1}{rgb}{0.99, 0.5, 0}
\definecolor{emmablue1}{rgb}{0, 0.25, 0.5}
\definecolor{emmaskyblue1}{rgb}{0.4, 0.8, 0.95}
\definecolor{emmapurple1}{rgb}{0.5, 0.25, 0.6}

\pgfplotscreateplotcyclelist{mycolorlist3}{%
emmagreen1\\%
emmaorange1\\%
emmablue1\\%
emmapurple1\\%
emmaskyblue1\\%
}

\pgfplotscreateplotcyclelist{mygraylist}{%
black\\%
gray1, dashed\\%
gray3\\%
gray2\\%
gray2,dashed\\%
gray3,dashed\\%
}

\renewcommand{\baselinestretch}{0.969}

\title{\vspace*{18pt} On Fundamental Limitations of Dynamic Feedback Control in Regular Large-Scale Networks }

\author{ {Emma Tegling, Partha Mitra, Henrik Sandberg and Bassam Bamieh}
 
 \thanks{ E. Tegling and H. Sandberg are with the School of Electrical Engineering and Computer Science and the ACCESS Linnaeus Center, KTH Royal Institute of Technology, SE-100 44 Stockholm, Sweden {(\tt tegling, hsan@kth.se)}. P. Mitra is with Cold Spring Harbor Laboratory, Cold Spring Harbor, NY, USA, 11724 {\tt (mitra@cshl.edu)}. B. Bamieh is with the Department of Mechanical Engineering at  University of California Santa Barbara, Santa Barbara, CA, USA, 93106. {\tt (bamieh@engr.ucsb.edu) } Corresponding author: E. Tegling.} 
  \thanks{This work was supported in part by the Swedish Research Council through grants 2013-5523 and 2016-00861, and by the National Science Foundation through the INSPIRE grant PHY-1344069 and EECS -1408442.}

}
        
\begin{document}
\maketitle
\thispagestyle{empty}
%\pagestyle{empty}

%%%%%%%%%%%%%%%%%%%%%%%%%%%%%%%%%%%%%%%%%%%%%%%%%%%%%%%%%%%%%%%%%%%%%%%%%%%%%%%%%%%
% Abstract
%%%%%%%%%%%%%%%%%%%%%%%%%%%%%%%%%%%%%%%%%%%%%%%%%%%%%%%%%%%%%%%%%%%%%%%%%%%%%%%%%%%

\begin{abstract}
We study fundamental performance limitations of distributed feedback control in large-scale networked dynamical systems. Specifically, we address the question of whether dynamic feedback controllers perform better than static (memoryless) ones when subject to locality constraints. We consider distributed {linear} consensus and vehicular formation control problems modeled over {toric} lattice networks. {For the resulting spatially invariant systems we} 
study the large-scale asymptotics (in network size) of global performance metrics that quantify the level of network coherence. With static feedback {from relative state measurements,} such metrics are known to scale unfavorably in lattices of low spatial dimensions, preventing, for example, a 1-dimensional string of vehicles to move like a rigid object. We show that the same limitations in general apply also to dynamic feedback control that is locally of first order. This means that the addition of one local state to the controller gives a similar asymptotic performance to the memoryless case. 
This holds unless the controller can access noiseless measurements of its local state with respect to an absolute reference frame, in which case the addition of controller memory may fundamentally improve performance. {In simulations of platoons with 20-200 vehicles we show that the performance limitations we derive manifest as unwanted accordion-like motions. Similar behaviors are to be expected in any network that is embeddable in a low-dimensional toric lattice, and the same fundamental limitations would apply.   } To derive our results, we present a general technical framework for the analysis of stability and performance of {spatially invariant systems in the limit of large networks}.

\end{abstract}

%%%%%%%%%%%%%%%%%%%%%%%%%%%%%%%%%%%%%%%%%%%%%%%%%%%%%%%%%%%%%%%%%%%%%%%%%%%%%%%%%%%
% Introduction
%%%%%%%%%%%%%%%%%%%%%%%%%%%%%%%%%%%%%%%%%%%%%%%%%%%%%%%%%%%%%%%%%%%%%%%%%%%%%%%%%%%
\section{Introduction}
A central problem in the control of networked systems is to understand and quantify how architectural constraints on the controller affect global performance. A prototypical scenario is when sensing and actuation is distributed across a network,  but the controller's  connectivity is limited and localized, with the architectural constraints being described by a graph structure. In such settings, meaningful performance metrics are typically global in character as they involve aggregates of quantities from across the entire network. Natural questions arise as to how these global performance metrics are impacted by architectural constraints; for example, how does increasing the size of sensing neighborhoods, or altering the topological connectivity of the controller's network affect the best achievable performance?
These types of questions arise in a wide range of applications such as 
vehicle formation control, transportation systems, sensor networks, and electric power systems.

The question motivating the present work is whether {\em dynamic} feedback controllers perform better than static (that is, memoryless) ones in large-scale networks. {For {regular} networks, 
we study comparatively the large-scale asymptotics (in network size) of global performance of localized static versus dynamic feedback control. Previous work~\cite{Bamieh2012} has shown the dimensionality dependence of localized {\em static} state-feedback control in  {such} networks, with lower dimensional lattices having worse asymptotic performance scalings than higher dimensional ones. 
{This implies for example that} vehicular platoon formations, {which resemble} one-dimensional lattices, exhibit {limitations} in terms of the feasibility of constructing a formation that moves like a rigid object. 
An important question arises as to whether the use of localized {\em dynamic} feedback control may alleviate these limitations, 
that is, whether the controller's additional memory may compensate for the lack of global sensing. 
%\vspace{-0.25mm}

This question of static versus dynamic feedback is a version of an old question in the area of decentralized control~\cite{vsiljak1978large}. It can be motivated by recalling the following important fact about state feedback control. 
In fully centralized optimal linear quadratic control (e.g. LQR or state-feedback $H_\infty$ control), {\em static state feedback} is optimal. In other words, there is no additional advantage in using dynamic or time-varying controllers over static gains when the full state is available for feedback. This is, however, no longer true when architectural constraints are imposed on the controller, such as diagonal or banded structures~\cite{wang1973stabilization,willems1989time,anderson1981time}. In our case, architectural constraints corresponding to bandedness are imposed through a {periodic} lattice-network structure, which motivates the use of dynamic feedback in search for the best-achievable performance.

Overall, a common theme in the area of networked systems involves designing control and interaction rules for a given multi-agent system and then showing that these rules lead to the desired performance in terms of stability or robustness. A difficulty with this approach is that when performance guarantees for the proposed interaction rules cannot be found, it is not evident whether this is due to a fundamental limitation, or simply lack of ingenuity on part of the control designer. A systematic approach to understanding fundamental limitations is instead to solve optimal control problems and check whether the best achievable performance is acceptable or fundamentally limited. In fully centralized control, several such criteria based on the plant's unstable poles and zeros are well known~\cite{skogestad2007multivariable}. With architectural constraints on the controller, however, the optimal control problem is, with few exceptions, non-convex. 
For example, while it is known that large-scale centralized optimal control problems have an inherent degree of locality~\cite{Bamieh2002,Motee2008}, optimal design of a controller with a {\em prescribed} degree of locality is non-convex. Important exceptions to this include the subclass of funnel-causal and quadratically-invariant problems (see e.g.~\cite{rotkowitz2006characterization,bamieh2005convex}, and a more recent approach~\cite{wang2017separable}). 

For large-scale networks, an emerging approach to understanding fundamental limitations while overcoming the non-convexity difficulty is  to derive asymptotic bounds on the best achievable performance~\cite{Bamieh2012,Lin2012,  Patterson2014, SiamiMotee2015,Grunberg2016,Barooah2009,herman2015nonzero,herman2017scaling}. 
In particular, the approach taken in~\cite{Bamieh2012} is to study performance of distributed static state-feedback controllers with locality constraints. While the corresponding optimization problem is non-convex for any finite system size,  informative performance bounds are derivable in the limit of a large system.  The aforementioned studies have focused on the dependence of the best-achievable performance bounds on node dynamics and network topology (e.g. lattices and fractals of various dimensions). However, work has thus far been limited to {distributed} static feedback control which does not alter the structure of local dynamics. Hence, the possible impact of controller dynamics on performance bounds, which we address in this paper, has remained an open question.

The dynamic feedback controllers we consider are modeled with first order dynamics and can share their state locally.
The resulting control laws can be seen as generalized distributed proportional-integral (PI) controllers. 
Subclasses of such controllers have been proposed in the literature for the elimination of stationary control errors that arise through static feedback, which has made them relevant for distributed frequency regulation of power networks~\cite{Andreasson2014,  SimpsonPorco2013, Zhao2015}. Here, we consider two classes of systems with respectively first and second order local {integrator} dynamics; consensus and vehicular formation control problems. {In the latter, we limit the analysis to symmetric feedback interactions, corresponding to an undirected network graph.} {The systems are modeled over {a class of regular networks, specifically} toric lattices, with a fixed number of neighbor interactions.} In line with related work~\cite{Bamieh2012,Lin2012, Patterson2014, SiamiMotee2015,Grunberg2016}, we characterize performance through 
nodal variance measures that capture the notion of network \emph{coherence}, or the rigidity of the network {lattice} formation.

Our main result shows that the fundamental limitations in terms of asymptotic performance scalings that apply to localized static feedback {in these types of networks} in general carry over also to dynamic feedback where the controllers are locally first-order. 
This means that while additional memory in the controller may offer other advantages in architecturally constrained control problems, it will not alleviate the unfavorable scaling of nodal variance in low-dimensional lattice networks. 
An important exception to this result applies in formation control problems if the controller can access \emph{absolute} measurements of local velocity with respect to a global reference frame. 
In this case, a carefully designed dynamic feedback controller can {theoretically} achieve \emph{bounded} variance
for any lattice network, thereby enabling a vehicle platoon to move like a rigid object. 

{As the problems we consider are modeled over toric (i.e., periodic) lattice networks, the resulting systems are \emph{spatially invariant}, a class of systems described in detail in~\cite{Bamieh2002}. The topological restriction is a consequence of the aim of the study; to characterize performance scalings in network size. This requires a possibility to grow the network while preserving certain topological properties, such as locality. Hypercubic lattices (including those with periodicity) is one of a few {regular} graph families with such topological invariance properties. Others include triangular lattices and fractals, see e.g.~\cite{Patterson2014}. }

{The periodicity of the lattices allows feedback protocols to be defined using multidimensional circulant operators. These enable a tractable spectral characterization through Fourier analysis. At large system sizes, however, the periodic boundary condition will have little or no effect on behaviors in the interior of the network. This intuitive reasoning can be attributed to exponential spatial decay rates of local perturbations~\cite{Bamieh2002}. This in particular implies that the lack of coherence that our results predict for one-dimensional ring-shaped lattices will also be observed in, for example, vehicular platoons without the periodic boundary condition. This is also demonstrated through a case study in Section~\ref{sec:examples}.}

{We will next lay out some technical preliminaries for our analysis, and then proceed to} set up the problems with static versus dynamic feedback control in Section~\ref{sec:setup}. In Section~\ref{sec:perfmeas}, the performance measure is defined through the variance of nodal state fluctuations. This variance is shown to correspond to a (scaled) system \hn norm. As one of the main contributions of {this} work, Section~\ref{sec:finiteinfinite} introduces a novel framework for evaluating such~\hn norms and their asymptotic scalings {in spatially invariant systems}. {This framework} allows for an analysis of large classes of dynamic feedback protocols whose \hn norm expressions are otherwise intractable. {It is also useful for} analyzing stability of these protocols, to which we devote Section~\ref{sec:stability}. {We show here that several control designs inevitably destabilize the system as the networks grow large, rendering those designs \emph{inadmissible}. The performance of admissible feedback protocols is then analyzed} in Section~\ref{sec:evaluation}, where our main result is derived.
{In Section~\ref{sec:examples} we discuss practical implications of our results and present a numerical simulation.} We end in Section~\ref{sec:discussion} with a discussion of our findings as well as some open problems.

%%%%%%%%%%%%%%%%%%%%%%%%%%%%%%%%%%%%%%%%%%%%%%%%%%%%%%%%%%%%%%%%%%%%%%%%%%%%%%%%%%%

\subsection{Preliminaries and notation}
\label{sec:preliminaries}
Throughout this paper, we will consider problems over the undirected $d$-dimensional torus $\mathbb{Z}_L^d$ with a total of $N = L^d$ nodes and $d$ assumed finite. In the one-dimensional case ($d = 1$), $\mathbb{Z}_L$ is simply the $L$ node ring graph, which we can represent by the set of integers $\{-\frac{L}{2},\ldots,0,1,\ldots,\frac{L}{2}-1 \}$ mod $L$ for $L$ even, and $\{-\frac{L-1}{2},\ldots,0,1,\ldots,\frac{L-1}{2}\}$ mod $L$ for $L$ odd. $\Znd$ is the direct product of $d$ such rings. It will also be useful to define the infinite $d$-dimensional lattice $\mathbb{Z}^d$, which is the direct product of $d$ copies of the integers. 

We define real-valued function arrays over this network, such as $a:~\mathbb{Z}_L^d \mapsto \mathbb{R}$, where we will use multi-index notation to denote the $k^{\mbox{th}}$ array entry $a_k = a_{(k_1,\ldots, k_d)}$. Similarly, we denote the state at node $k= (k_1,\ldots, k_d)$ in the $d$-dimensional torus as $x_{(k_1,\ldots,k_d)}(t),$
which is a scalar in $\mathbb{C}$ in the consensus problems and a vector-valued signal in $\mathbb{R}^d$ in the vehicular formation problems. We will in most cases omit the time dependence in the notation.

Linear operators, denoted by upper case letters, will be used to define multi-dimensional circular convolutions with function arrays over  $\mathbb{Z}_L^d $. For example, the convolution operator~$A$ associated with the array~$a$ is defined as follows:
\begin{equation}\label{eq:convdef}
h \! = \!Ax ~ \Leftrightarrow ~ h_{(k_1,\ldots,k_d)} \! = \!\! \!\! \sum_{(l_1,\ldots,l_d) \in \mathbb{Z}_L^d} \!\! \!\! a_{(k_1,\ldots,k_d) - (l_1,\ldots,l_d)} x_{(l_1,\ldots,l_d)},
\end{equation}
or, in short, $h_k = (Ax)_k = \sum_{l \in \mathbb{Z}_L^d } a_{k-l} x_l.$

In cases where the state $x \in \mathbb{R}^d$, the array element $a_k$ is a $d \times d$ matrix, which in this paper will be assumed to be diagonal due to coordinate decoupling. 
The addition of multi-indices in the $\Znd$ arithmetic is done as follows: $k + l = (k_1,\ldots,k_d) + (l_1,\ldots,l_d)  = (k_1+l_1,\ldots,k_d+l_d)_{\mathrm{mod}L}$. Here, mod $L$ implies that the operation is circulant. Note that all feedback operators considered in this paper are \emph{spatially invariant} with respect to~$\Znd$, and can therefore be represented by convolution operators with single-index arrays as in \eqref{eq:convdef}.

The spatial discrete Fourier transform (DFT) of the array~$a$ will be denoted with $\hat{a}$, and we will use the letter~$n$ to denote the index, or \emph{wavenumber}, of the spatial Fourier transform. For example, the function array $a_{(k_1,\ldots, k_d)}$ has $\hat{a}_{(n_1,\ldots, n_d)}$ as its Fourier transform, where the wavenumber $(n_1,\ldots, n_d)$ can be thought of as a spatial frequency variable. Throughout this paper, we will use the DFT that is defined as: %for function arrays $a:~\mathbb{Z}_L^d \mapsto \mathbb{C}$ as:
\begin{equation}
\label{eq:DFTdef}
\hat{a}_n := \sum_{k \in \Znd} a_k e^{-j\frac{2\pi}{L}n \cdot k},
\end{equation}
where $j = \sqrt{(-1)}$ denotes the imaginary number and $n \cdot k = n_1k_1 + \cdots + n_dk_d$.  

Function arrays can also be defined over the infinite $d$-dimensional lattice $\mathbb{Z}^d$. We then use the subscript $\infty$ for the array, as in $a_\infty$, with entries $a_{(k_1,\ldots, k_d)}$ for $k \in \mathbb{Z}^d$. The corresponding convolution operator is denoted $A_\infty$. The $Z$-transform of $a_\infty$ evaluated on the unit circle $e^{j\theta}$ is:
\begin{equation}
\label{eq:ZtransformPrel}
\hat{a}_\infty(\theta)  := \sum_{k \in \mathbb{Z}^d} a_k e^{-j \theta \cdot k},
\end{equation} 
where $\theta = (\theta_1, \ldots, \theta_d)$ denotes a spatial frequency, which takes values in the multivariable rectangle $\mathcal{R}^d:=[-\pi,\pi]^d$. 

We will use the term \emph{(generalized) Fourier symbol} of convolution operators for the DFT or $Z$-transform of the corresponding function array. For example, $\hat{a}$ in \eqref{eq:DFTdef} is the Fourier symbol of the operator $A$. The values that $\hat{a}$ takes are exactly the eigenvalues of $A$. In cases where $a$ is matrix valued, the eigenvalues of $A$ are the union of all eigenvalues of $\hat{a}_{(n_1,\ldots, n_d)}$ as $(n_1,\ldots, n_d)$ runs through $\mathbb{Z}_L^d$.

 In this paper, we derive what we call \emph{asymptotic scalings } of certain performance measures with respect to the network size. The symbol~$\sim$ is used throughout to denote scalings in the following manner: \begin{equation}
 \label{eq:scalingdef}
 u(N) \sim v(N)~~\Leftrightarrow
 ~~\underline{c}v(N) \le u(N) \le \bar{c}v(N),
 \end{equation}
 for any $N$, where the fixed constants $\underline{c}, \bar{c}$ are independent of the variable $N$. When a scaling is said to hold \emph{asymptotically}, the relation~\eqref{eq:scalingdef} holds for all $N> \bar{N}$ for some fixed $\bar{N}$.

%%%%%%%%%%%%%%%%%%%%%%%%%%%%%%%%%%%%%%%%%%%%%%%%%%%%%%%%%%%%%%%%%%%%%%%%%%%%%%%%%%%
% Problem setup
%%%%%%%%%%%%%%%%%%%%%%%%%%%%%%%%%%%%%%%%%%%%%%%%%%%%%%%%%%%%%%%%%%%%%%%%%%%%%%%%%%%
\section{Problem Setup}
\label{sec:setup}
We now formulate models for two types of problems: \emph{consensus} and \emph{vehicular formations}. In the consensus problem there is a local, scalar information state at each network site, while there are two such states (position and velocity) in the vehicular formation case. 
For both models,  we introduce a \textit{static controller}, as considered in \cite{Bamieh2012}, which we will compare to a \textit{dynamic controller} with an auxiliary memory state at each network site, see Fig.~\ref{fig:structure}. 

\begin{figure}
\centering
\includegraphics[width = 0.42\textwidth]{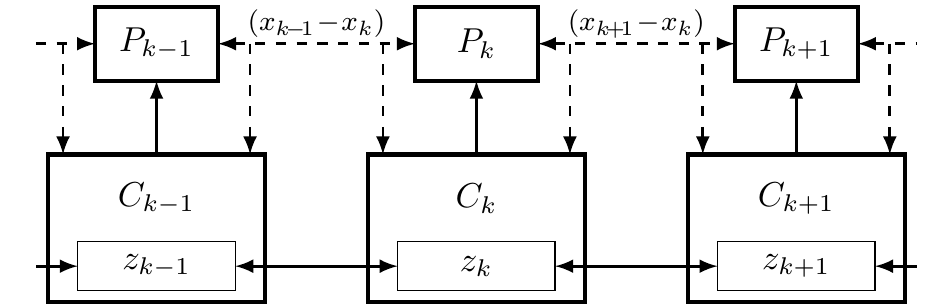}
%Tikz original is called controllerfig_new2
\caption{Structure of the controller $\{C_k\}$ and plant $\{P_k\}$ interactions. For diagrammatic simplicity only nearest neighbor interaction is depicted, though our analysis is applicable to any fixed number of neighbor interactions.  Dashed arrows indicate relative state measurements and interactions. The controller states $\{z_k\}$, rather than just their relative values, can be shared between sub-controllers.   }
\label{fig:structure}
\end{figure}

\subsection{Consensus} 
\label{sec:cons_setup}
We first consider the first-order consensus algorithm in continuous time over the discrete torus $\mathbb{Z}_L^d$. The single-integrator dynamics at each site $k$ in the network is given by 
\begin{equation}
\label{eq:consensus1}
\dot{x}_k = u_k + w_k, ~~ k \in \mathbb{Z}_L^d,
\end{equation}
where $u_k$ denotes the control signal. The process disturbance~$w_k$, modeling random insertions and deletions, is mutually uncorrelated across nodes. Throughout this paper, we model this disturbance as zero mean white noise\footnote{We refer to ``white noise'' in continuous time as a stationary zero-mean stochastic process with autocorrelation $\mathbb{E}\{w(\tau)w^*(t)\} = \delta(t-\tau)I$, where $\delta(t)$ denotes the Dirac delta distribution. This idealized process can be thought of as the time derivative of a Brownian motion, $dB/dt$, although such a derivative does not formally exist, see \cite[Theorem 4.1]{AstromBook}.}. 

We now introduce the two types of linear time-invariant feedback control for the system \eqref{eq:consensus1}.

\subsubsection{Static feedback}
In the case of static feedback, the control input is a linear function of the current network state: %signal depends directly on the current state, such that
\begin{equation}
u_k = (Fx)_k.
\label{eq:consensusstaticcontrol}
\end{equation}
The feedback operator $F$, %(with the associated array $f$ defined on $\Znd$) 
can be suitably designed to fulfill the control objectives. A common example of such a control scheme is the one where the control signal at each node is the weighted average of the differences between that node and its 2$d$ neighbors, that is,
\begin{align} \nonumber
 u_k = &\tilde{f} [ (x_{(k_1-1,\ldots,k_d)} - x_k) +  (x_{(k_1+1,\ldots,k_d)} - x_k) + \cdots\\
   &+(x_{(k_1,\ldots,k_d-1)} - x_k) +  (x_{(k_1,\ldots,k_d+1)} - x_k) ], \label{eq:standardconsensus}
\end{align}
where $\tilde{f}$ is a positive scalar.
The algorithm \eqref{eq:standardconsensus} will be referred to as the \emph{standard consensus algorithm}. The associated function array is:
\begin{equation}
\label{eq:standardarray}
f_{(k_1,\ldots,k_d)} = \begin{cases} -2d\tilde{f} & k_1 = \cdots = k_d = 0 \\ \tilde{f} & k_i = \pm 1,~\mathrm{and}~k_j = 0, ~\mathrm{for}~ i \neq j \\ 0 & \mathrm{otherwise} .\end{cases}
\end{equation}

In the general case, we can write the consensus algorithm~\eqref{eq:consensus1} with static feedback as
\begin{equation}
\label{eq:consensusstatic}
\dot{x} = Fx + w.
\end{equation}

\subsubsection{Dynamic feedback}
To model dynamic feedback, we let the controller have access to an auxiliary controller state~$z_{(k_1,\ldots,k_d)}$, which is a scalar at each network site $k$:
\begin{align*}
u_k &= z_k + (Fx)_k \\
\dot{z}_k & = (Az)_k +(Bx)_k,
\end{align*}
where $A,~B,~F$ are linear feedback operators, the properties of which will be discussed shortly.  We can now write the consensus algorithm \eqref{eq:consensus1} with dynamic feedback as:
\begin{align}
\begin{bmatrix}
\dot{z} \\ \dot{x}
\end{bmatrix} = \begin{bmatrix}
A & B\\ I& F
\end{bmatrix} \begin{bmatrix}
z \\ x
\end{bmatrix} + \begin{bmatrix}
0 \\ I
\end{bmatrix}w  .
\label{eq:consensusdynamic}
\end{align}

\subsection{Structural assumptions for the consensus problem}
\label{sec:cons_ass}

We now list the assumptions imposed on the system and on the feedback operators $A,~B, ~F$ in the consensus algorithm. Assumptions \ref{ass:spatialinvariance}-\ref{ass:locality} will also carry over to the vehicular formation problems. 
\begin{ass}[Spatial invariance]
\label{ass:spatialinvariance}
All feedback operators are spatially invariant and fixed with respect to $\Znd$, and are therefore circular convolution operators, as defined in \eqref{eq:convdef}.

For example, the standard consensus algorithm \eqref{eq:standardconsensus} on the 1-D ring graph $\mathbb{Z}_L$ can be written as the convolution of the state $x$ with the array $f=\{0, \ldots,0,\tilde{f}, -2\tilde{f} ,\tilde{f} ,0,\ldots,0 \}$.

\end{ass}
\begin{ass}[Locality]
\label{ass:locality}
All feedback operators use only local information from a neighborhood of width $2q$, where $q$ is independent of $L$. For the function array~$f$ associated with the operator~$F$, this means that
\begin{equation}
f_{(k_1,\ldots,k_d)} = 0~~~\mbox{if}~~|k_i|>q,
\label{eq:localitycondition}
\end{equation}
for any $i \in \{1,2,\ldots,d\}$. The same condition holds for all other operators. The situation is illustrated in  Fig.~\ref{fig:arrayfig}.
\end{ass}
\begin{ass}[Relative state measurements]
\label{ass:relative}
%We assume that 
All controllers can only access \emph{relative} measurements of the physical state $x$. Hence, the feedback can only involve differences between states of neighboring nodes. This means that each term of the form $\tilde{f}x_k$ in the convolution $Fx$ is accompanied by another term $-\tilde{f}x_l$, for some other index $l$, so that we obtain $\tilde{f}(x_k - x_l)$.
In particular, this implies that the operators~$F$ and $B$ in~\eqref{eq:consensusstatic} and~\eqref{eq:consensusdynamic} have the property
\begin{equation}
\label{eq:relativedef}
\sum_{k \in \Znd} f_k = 0, ~~\sum_{k \in \Znd} b_k = 0.
\end{equation} 
Since the state $z$ is internal to each controller, %and can therefore always be accessed, 
we need \textit{not} impose this requirement on $A$ in~\eqref{eq:consensusdynamic}. 
\end{ass}

\begin{figure}[t!]
\centering
\includegraphics[width = 0.40\textwidth]{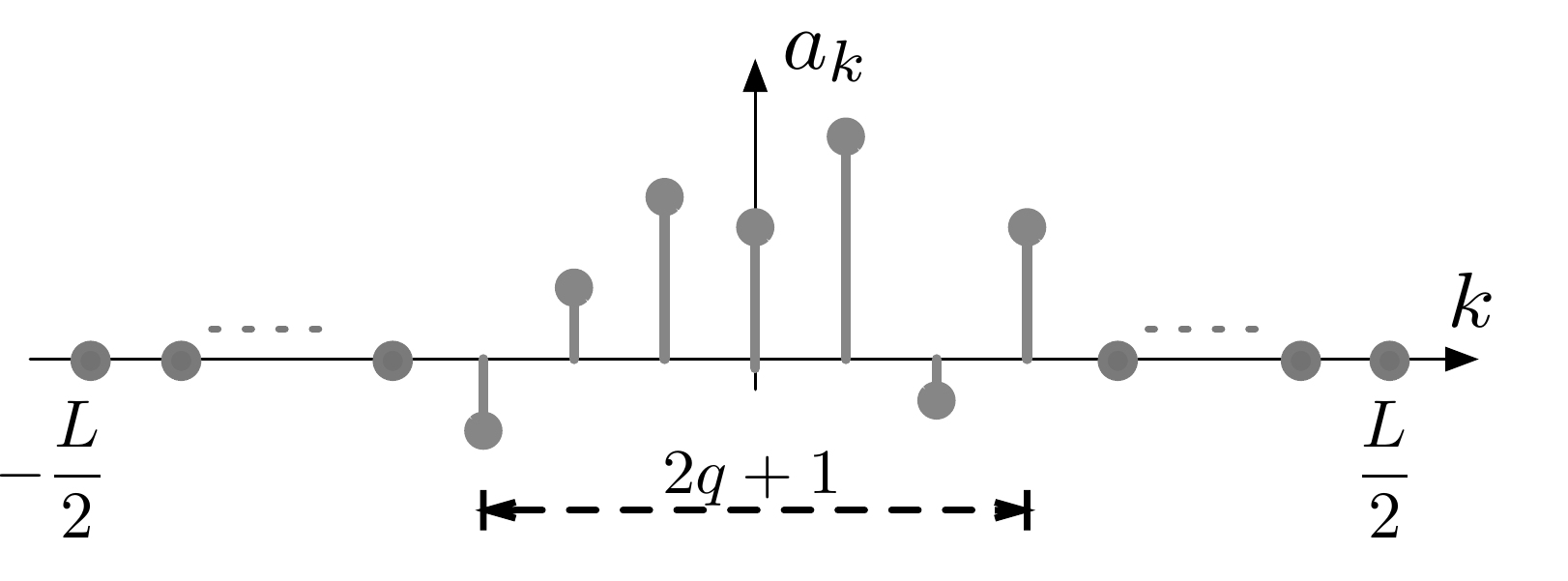}
\caption{Spatial interactions are defined by convolution with an array $\{a_k\}$. As the lattice size $L$ increases, the locality property ($a_k=0$ for $|k|>q$) insures that the interactions 
are unambiguously defined.}% as $L\rightarrow\infty$.}  
\label{fig:arrayfig}
\end{figure}

\subsection{Vehicular formations}
\label{sec:veh_setup}
For the vehicular formation problem, consider $N = L^d$ identical vehicles arranged in the $d$-dimensional torus $\Znd$. The double integrator dynamics at each site $k = (k_1, \ldots, k_d)$ in the torus is then
\begin{equation}
\ddot{x}_k = u_k + w_k,
\label{eq:veh1}
\end{equation}
where, as before, $u_k$ is the control signal and $w_k$ is white process noise, which models random forcings at each site. 

The position vector $x_k = [x_k^1 ~ \cdots ~ x_k^d]^T$ at each network site, and its time derivative, the velocity vector $v_k= [v_k^1 ~ \cdots ~ v_k^d]^T$, are both $d$-dimensional vectors. Without loss of generality, we will assume that they each represent absolute deviations from a desired trajectory $\bar{x}_k$ and constant heading velocity $\bar{v}$, with \[\bar{x}_k := \bar{v}t + k\Delta_x.\] Here, $\Delta_x$ is the constant spacing between the vehicles in $\Znd$. 

In analogy to the consensus case, we now introduce the two types of linear feedback control for the system \eqref{eq:veh1}. 

\subsubsection{Static feedback}
The control input is here assumed to be full state feedback that is linear in the variables $x$ and $v$:
\[u_k = (Fx)_k + (Gv)_k. \]
An example of such feedback is the combined look-ahead and look-behind controller in a 1-D string:
%\begin{align} \nonumber
%u_k =& f_+(x_{k+1} - x_k) + f_-(x_{k-1}-x_k) + g_+(v_{k+1} - v_k) +\\ \label{eq:vehfeedbackexample}
% & +g_-(v_{k-1} - v_k) - f_ox_k - g_ov_k,
%\end{align}
\begin{align} \nonumber
u_k =& f_+(x_{k+1} - x_k) + f_-(x_{k-1}-x_k) + g_+(v_{k+1} - v_k) +\\ \label{eq:vehfeedbackexample}
 & +g_-(v_{k-1} - v_k) - g_ov_k,
\end{align}
where the $g$'s and $f$'s are positive design parameters. {If $g_o$ is zero, this control law satisfies Assumption \ref{ass:relative} of relative state measurements. If $g_o \neq 0$, we will refer to that term in the feedback law as \textit{absolute feedback} from velocity.} 

{In practice, absolute velocity measurements can be made available through a speedometer. The presence of viscous damping can also be treated as a special case of absolute velocity feedback. The model \eqref{eq:veh1} can then be modified so that $\dot{v}_k = -\mu v_k + u_k + w_k,$
%\begin{equation}
%\label{eq:viscous}
%\ddot{x}_k = -\mu \dot{x}_k + u_k + w_k,
%\end{equation}
where $\mu\ge 0$ is the drag coefficient.  Comparing this to \eqref{eq:vehfeedbackexample} we can identify $\mu$ with $g_o$.} 

{We will not consider the case where absolute feedback is available from the position~$x_k$ but not from the velocity~$v_k$. Such a scenario would correspond to vehicles accessing absolute position measurements via e.g. GPS, yet lacking the ability to derive their absolute velocity from those measurements. See also Remark~\ref{rem:absposition}. }

%Suppressing the spatial index of all variables, 
In summary, the vehicular formation algorithm \eqref{eq:veh1} with the static feedback law becomes
\begin{equation}
\label{eq:vehiclestatic}
\begin{bmatrix}
\dot{x} \\ \dot{v}
\end{bmatrix} = \begin{bmatrix}
0 & I\\ F& G
\end{bmatrix} \begin{bmatrix}
x \\ v
\end{bmatrix} + \begin{bmatrix}
0 \\ I
\end{bmatrix}w.  
\end{equation}

\subsubsection{Dynamic feedback}
To model the dynamic feedback laws, we introduce the auxiliary controller state $z_k$ at each network site $k$, which is a $d$-dimensional vector containing a memory of past position and velocity errors in each coordinate direction. We get: 
\begin{align*}
u_k &= z_k +  (Fx)_k + (Gv)_k \\
\dot{z}_k &= (Az)_k +(Bx)_k + (Cv)_k.
\end{align*}
An example of dynamic feedback control for double integrator systems is \textit{distributed-averaging proportional-integral }(DAPI) control, which has received much recent attention in the context of coupled oscillator systems and control of so-called microgrids~\cite{Andreasson2014, SimpsonPorco2013, Zhao2015}. Such systems are analogous to the vehicular formation problem under certain assumptions, such as absolute velocity feedback. One DAPI control algorithm is:
\begin{equation}\label{eq:DAPIexample}
\begin{aligned}
u_k = &z_k + f_+(x_{k+1} - x_k) + f_-(x_{k-1}-x_k)  - g_ov_k\\ 
\dot{z}_k = &a_+(z_{k+1} - z_{k}) + a_-(z_{k-1} - z_k) - c_ov_k
\end{aligned}
\end{equation}
\iffalse
\begin{align} \nonumber
u_k = &z_k + f_+(x_{k+1} - x_k) + f_-(x_{k-1}-x_k)  - g_ov_k\\ \label{eq:DAPIexample}
\dot{z}_k = &a_+(z_{k+1} - z_{k}) + a_-(z_{k-1} - z_k) - c_ov_k
\end{align}\fi
where the operator $A$ achieves a weighted averaging of the internal state $z$ across nodes, which prevents unfavorable drift in the memory states at different nodes. Such drift would, in practice, de-stabilize the system if $A = 0$, in which case~\eqref{eq:DAPIexample} reduces to a \textit{decentralized} proportional-integral~(PI) controller with respect to the velocity, see e.g.~\cite{Andreasson2014}. 

In general, we can write the equations of motion for the closed loop system with dynamic feedback as:
\begin{align}
\label{eq:vehicledynamic}
\begin{bmatrix} 
\dot{z} \\ \dot{x}  \\ \dot{v}
\end{bmatrix} &= \begin{bmatrix}
A & B & C \\ 0& 0 & I \\ I&F&G
\end{bmatrix} \begin{bmatrix} 
z \\ x \\v
\end{bmatrix} + \begin{bmatrix}
0 \\ 0 \\ I
\end{bmatrix}w.
\end{align}
%where we have again suppressed the spatial indices for all variables. 

\subsection{Structural assumptions for the vehicular formation problem}
\label{sec:veh_ass}
For the vehicular formation systems, we impose the following assumptions\emph{ in addition to} Assumptions~\ref{ass:spatialinvariance}-\ref{ass:locality} {above.}%listed in Section~\ref{sec:cons_ass}:

\begin{ass}[Relative position measurements] 
\label{ass:relabs} {The controllers can only access relative measurements of the position states~$x$. This means that the operators $F$ and $B$ in~\eqref{eq:vehiclestatic} and~\eqref{eq:vehicledynamic} have the property~\eqref{eq:relativedef}.}
\end{ass}

\begin{rem}
{We will both consider the case where the velocity feedback operators~$G$ and $C$ have the relative measurement property~\eqref{eq:relativedef} as well as when they do not. We refer to these cases as, respectively, \emph{relative} and \emph{absolute velocity feedback}.}
\end{rem}

\begin{ass}[Reflection symmetry]
\label{ass:symmetry}
The interactions between the vehicles on $\Znd$ are symmetric around each site~$k$. This implies that the arrays associated with the operators $A,B,C,F,G$ have even symmetry, so that for each array element $ f_{(k_1,\ldots,k_d)}= f_{(-k_1, \ldots,-k_d)}$.
% for each term like $\alpha f_{(k_1,\ldots,k_d)}$ there is a term $\alpha f_{(-k_1, \ldots,-k_d)}$. 
For example, in~\eqref{eq:DAPIexample} this condition requires $a_+ = a_-$, $b_+ = b_-$, $f_+ = f_-$ and $g_+ = g_-$. 

A particular implication of this assumption is that the Fourier symbols of the operators will be real valued. 

The property of reflection symmetry will be relevant (but not enforced) also in the consensus case. By slight abuse of terminology, we will in the following refer to a feedback operator as \emph{symmetric} if the associated array has this property, and \emph{asymmetric} if it does not. 
\end{ass}

\begin{ass}[Coordinate decoupling]
\label{ass:decoupling}
The feedback in each of the $d$ coordinate directions is entirely decoupled from the vector components in the other coordinates. Furthermore, the array elements associated with the operators $A,B,C,F,G$ are isotropic. %identical in each of the coordinate directions. 
By this assumption, the array elements $a,b,c,f,g$ are diagonal and the convolution in~\eqref{eq:convdef} will turn into $d$~decoupled, identical, scalar convolutions. 
\end{ass}
{While Assumptions~\ref{ass:spatialinvariance}--\ref{ass:relabs} are important for the upcoming analysis, Assumptions~\ref{ass:symmetry}--\ref{ass:decoupling} are mainly made to simplify the  calculations.}

%%%%%%%%%%%%%%%%%%%%%%%%%%%%%%%%%%%%%%%%%%%%%%%%%%%%%%%%%%%%%%%%%%%%%%%%%%%%%%%%%%%
% Perf meas & main result
%%%%%%%%%%%%%%%%%%%%%%%%%%%%%%%%%%%%%%%%%%%%%%%%%%%%%%%%%%%%%%%%%%%%%%%%%%%%%%%%%%%

\section{Performance Measure and Main Result}
\label{sec:perfmeas}
In this paper, we are concerned with the performance of the consensus and vehicular formation problems in terms of the amount of global ``disorder'' of the system at steady state. This can be quantified as the steady state variance of nodal state fluctuations, which are caused by persistent stochastic disturbances. In particular, we are interested in the \emph{scaling} of this performance measure with the system size, as it grows asymptotically. We call a system that exhibits a better scaling \emph{more coherent} than a system with bad scaling, as the former will form a more rigid formation when the system grows. If the scaling is such that the variance per node is bounded, the system is said to be \emph{fully coherent}.

We adopt the approach in \cite{Bamieh2012} to define the relevant performance measure. %, but restrict our attention to one of the three performance measures evaluated there. 
Consider first a general linear MIMO system driven by zero mean white noise $w$ with unit intensity:% (i.e., $\mathbb{E}\{w(\tau)w^*(t)\} = \delta(t-\tau)I$):
\begin{subequations}
\label{eq:generalSS}
\begin{align} \label{eq:generalSSsystem}
\dot{\psi}(t) &= \mathcal{A}\psi(t)+\mathcal{B}w(t)\\  \label{eq:generalSSoutput}
y(t) & = \mathcal{C}\psi(t).
\end{align}
\end{subequations}
In our case, equation \eqref{eq:generalSSsystem} represents, for example, the feedback system in \eqref{eq:vehicledynamic}, for which a performance output $y(t)$ as in~\eqref{eq:generalSSoutput} will be defined shortly.

Provided that the system \eqref{eq:generalSS} is input-output stable, its squared \hn norm from $w$ to $y$ is finite and can be interpreted as the total steady state variance of the output, that is,
\begin{equation}
\label{eq:variancedef}
\mathbf{V}_N: = \sum_{k \in \Znd} \lim_{t\rightarrow \infty} \mathbb{E} \{y_k^*(t) y_k(t)\}.
\end{equation}
Throughout this paper, we are considering spatially invariant systems over the discrete torus $\Znd$. This implies that the output variance $\mathbb{E}\{y_k^*(t)y_k(t)\}$ will be independent of the site~$k$.
We obtain this {steady state} \emph{per-site variance }by simply dividing the total \hn norm by the system size~$N = L^d$:
\begin{equation}
\label{eq:pernodevariance}
\Vn = \lim_{t\rightarrow \infty}\mathbb{E}\{y_k^*(t)y_k(t)\} = \frac{\mathbf{V}_N}{N}.
\end{equation}

We next define the relevant output measurement that will be used throughout the paper:
%We next define the considered performance measure:
\begin{definition}[Deviation from average performance measure]
%The deviation of each state from the average of all states is measured as:
\begin{equation}
\label{eq:performancedef}
y_k := x_k - \frac{1}{N} \sum_{l \in \Znd} x_l.
\end{equation}
In operator form, this becomes 
\begin{equation}
\label{eq:Hdavdef}
y = (I - \frac{1}{N}J_{\mathbf{1}})x =: Hx, 
\end{equation} where $J_{\mathbf{1}}$ is the convolution operator corresponding to the array with all elements equal to $1$.
\end{definition} 

\begin{rem}
\label{rem:BIBO}
It is well known that the consensus type dynamics considered in this paper typically have a single marginally stable mode at the origin corresponding to the motion of the average (this is a consequence of Assumption~\ref{ass:relative} of relative measurements). The \hn norm~\eqref{eq:variancedef} is only finite if this mode is unobservable from the system output. Here, the output operator~$H$ has the relative measurement property~\eqref{eq:relativedef}, that is, $\sum_{k \in \Znd}h_k = 0$, implying that the average mode is indeed unobservable. Provided remaining system modes are stable, $\mathbf{V}_N$ in~\eqref{eq:variancedef} will thus be finite for any finite system size~$N$; a condition equivalent to bounded-input, bounded-output (BIBO) stability.
\end{rem}

\subsection{Performance scalings with static and dynamic feeedback}
The main objective of this paper is to determine whether dynamic feedback %on the forms \eqref{eq:consensusdynamic} and \eqref{eq:vehicledynamic} 
may improve performance compared to the static feedback laws %\eqref{eq:consensusstatic} and \eqref{eq:vehiclestatic} 
that were also evaluated in \cite{Bamieh2012}. {The} following sections {will} introduce the methodology that {is used} to establish asymptotic scalings of {performance}. At this point, we summarize our main results as follows.

\begin{theorem}[Asymptotic performance scalings]
\label{thm:mainresult}
{ Consider the consensus problem %defined in Section~\ref{sec:cons_setup} 
under Assumptions~\ref{ass:spatialinvariance}--\ref{ass:relative} and the vehicular formation problem %defined in Section~\ref{sec:veh_setup} 
under Assumptions~\ref{ass:spatialinvariance}, \ref{ass:locality}, and~\ref{ass:relabs}--\ref{ass:decoupling}. The steady state per-site variance~$V_N$ defined in~\eqref{eq:pernodevariance} then scales asymptotically as follows:}
\begin{enumerate}
\item Consensus 
%\begin{enumerate}
%\item 

\emph{Static} feedback \emph{or} \emph{dynamic} feedback 
\begin{equation}
\label{eq:resultscalingcons}
\Vn \sim \frac{1}{\beta}\begin{cases} N & d = 1\\  \mbox{log}N & d = 2 \\  1& d \ge 3, \end{cases}  
\end{equation}
%\end{enumerate}
\item Vehicular formations
\begin{enumerate}
%\item \emph{Static} feedback \emph{or} \emph{dynamic} feedback{: relative feedback}
\item {Relative feedback:} 

\emph{Static} feedback \emph{or} \emph{dynamic} feedback
\begin{equation}
\label{eq:resultscalingvehbad}
\Vn \sim \frac{1}{\beta^2}\begin{cases} N^3 & d = 1\\  N & d = 2 \\  N^{1/3}& d = 3 \\  \mbox{log}N & d = 4 \\ 1 & d\ge 5, \end{cases} 
\end{equation}
% \item \emph{Static} feedback{: absolute velocity feedback} {(but relative position feedback)}
\item {Absolute velocity (but relative position) feedback:}% {(relative position feedback): } 

\emph{Static} feedback
\begin{equation}
\label{eq:resultscalingvehbetter}
\Vn \sim \frac{1}{\beta}\begin{cases} N & d = 1\\  \mbox{log}N & d = 2 \\  1& d \ge 3, \end{cases}
\end{equation}
\emph{Dynamic} feedback
%\item  \emph{Dynamic} feedback{: absolute velocity feedback} {(but relative position feedback)}
%
%\noindent Static feedback \emph{or} dynamic feedback with \emph{absolute} position and \emph{absolute} velocity feedback
\begin{equation}
\label{eq:resultscalingvehbest}
\Vn \sim 1,  
\end{equation}
\end{enumerate}
\end{enumerate}
where $N = L^d$ is the network size, {$\beta = \max\{||f||_\infty,||g||_\infty\}$} is an algorithm parameter reflecting {the magnitude of feedback gains}, and the {symbol $\sim$ denotes} scaling up to a factor that is independent of $N$ and $\beta$~{in the manner defined in~\eqref{eq:scalingdef} }.
\end{theorem}
\vspace{0.5mm}

Therefore, if only relative state measurements are available (Assumption~\ref{ass:relative}), no dynamic feedback laws on the forms~\eqref{eq:consensusdynamic} and~\eqref{eq:vehicledynamic} exhibit better coherence properties than static, memoryless feedback under the given assumptions.
%In the vehicular formation case, the same conclusion holds also in the case of absolute position and relative velocity feedback.

However, a dynamic feedback law can theoretically achieve full coherence in any spatial dimension using absolute feedback from velocities, even though position measurements are relative. {As previously shown in~\cite{Bamieh2012}, and as will be evident from the developments in Section~\ref{sec:evaluation},} a static feedback protocol {would} require absolute measurements of \emph{both} states to achieve the same performance. 

\begin{rem}
\label{rem:absposition}
{ The case with absolute position but relative velocity feedback is not considered here. %As already discussed, it is of less practical relevance than the converse.
The interested reader is referred to~\cite{Tegling2016Lic} where it is shown that $V_N$ then scales as in~\eqref{eq:resultscalingvehbetter} for both the static and dynamic feedback laws modeled in this paper. 
%static feedback and dynamic feedback on the form~\eqref{eq:vehicledynamic} both lead to the scaling of~$V_N$ in~\eqref{eq:resultscalingvehbetter}.  
In~\cite{Tegling2017}, however, an alternative controller with derivative action is also designed for this particular case, which is shown to give full coherence.  }
\end{rem}

%%%%%%%%%%%%%%%%%%%%%%%%%%%%%%%%%%%%%%%%%%%%%%%%%%%%%%%%%%%%%%%%%%%%%%%%%%%%%%%%%%%
% Limits
%%%%%%%%%%%%%%%%%%%%%%%%%%%%%%%%%%%%%%%%%%%%%%%%%%%%%%%%%%%%%%%%%%%%%%%%%%%%%%%%%%%

\section{The \hnnorm Density and Asymptotic Scalings}
\label{sec:finiteinfinite}
We now introduce the technical framework that will be used to determine the \hn performance scalings in Theorem~\ref{thm:mainresult}. This novel framework is based on the idea of mapping the {operators that define the }system dynamics onto an infinite lattice. Usually, \hn norms are calculated using traces of system Gramians that lead to sums involving system eigenvalues. In the limit of large systems, they can instead be estimated through integrals over a continuous function which we call the \hnnorm density. We show that simple properties of this \hnnorm density determine the asymptotic performance scalings.% of the \hn performance. 

\subsection{The limit from finite to infinite lattices}
%\subsection{Spatial convolution operators on infinite lattices}
\label{sec:operatorsinfinite}
All feedback operators considered in this paper define convolutions with \emph{local} arrays on $\Znd$, by Assumption~\ref{ass:locality}. Hence, for a given operator $A$ we have that $a_k = 0$ if $|k|>q$ for some fixed $q$. This means that any such array $a$ can be unambiguously re-defined on $\mathbb{Z}_{L'}^d$ for \emph{any} given $L'>2q$ by filling it with zero components wherever $|k|>q$. This also means that the same array can be used to define a convolution over the \textit{infinite} lattice $\mathbb{Z}^d$. As we shall see, such a re-definition proves useful when analyzing the systems asymptotically.
%think of the sites in the finite torus $\Znd$ as a subset of the points in the infinite lattice $\mathbb{Z}^d$.

Let $a$ be a local array defined over $\Znd$ and $a_\infty$ its counterpart defined on $\mathbb{Z}^d$, in which the elements $\{a_k\}$ have been filled out with zeros for $|k| > q$ up until infinity. The discrete Fourier transform (DFT) of $a$, denoted $\hat{a}_n$ is given by~\eqref{eq:DFTdef} in Section~\ref{sec:preliminaries}, while the $Z$-transform of the array $a_\infty$, denoted $\hat{a}_\infty(\theta)$, is given by~\eqref{eq:ZtransformPrel}. 

Comparing \eqref{eq:DFTdef} with \eqref{eq:ZtransformPrel} it is clear that the DFT of $a$ is simply 
%the samples of the $Z$-transform of $a_\infty$ at the grid points $\theta = \frac{2\pi}{L}n$, for $n \in \Znd$. That is,
sub-samples of the $Z$-transform of~$a_\infty$:
\begin{equation}
\label{eq:sampling}
\hat{a}_n = \hat{a}_\infty\left(\frac{2\pi}{L} n \right), ~~n \in \Znd.
\end{equation}
Given that we are interested in system behaviors as $N \rightarrow \infty$, it will be convenient to consider these $Z$-transforms of operators over the infinite lattice $\mathbb{Z}^d$, and their behavior in the continuous spatial frequency variable $\theta \in \mathcal{R}^d$, rather than the DFTs at discrete spatial wavenumbers. 

For this purpose, let us take the general state space system~\eqref{eq:generalSS} and map the system operators $\mathcal{A}, ~\mathcal{B}, ~\mathcal{C}$ onto $\mathbb{Z}^d$ to obtain $ \mathcal{A}_\infty, \mathcal{B}_\infty, \mathcal{C}_\infty $. For example, in the system \eqref{eq:consensusstatic}, we have $\mathcal{A} = F$. If we let $F$ represent the standard consensus algorithm \eqref{eq:standardconsensus}, then $\mathcal{A}_\infty = F_\infty$ has the associated function array $f_\infty$, defined just as in \eqref{eq:standardarray}, but filled with infinitely many zero components for $|k_i|> 1$, $i = 1,\ldots,d$. %\emma{Put the preceding example in an ``Example-environment''?}

By virtue of the spatial invariance property, $ \mathcal{A}_\infty, ~\mathcal{B}_\infty$ and $\mathcal{C}_\infty $ are circulant convolution operators and the $Z$-transform can be used to \emph{(block)} \emph{diagonalize} them, see \cite{Bamieh2002}. Then, at each $\theta \in \mathcal{R}^d  $, we obtain the matrix-valued transforms $ \hat{\mathcal{A}}_\infty(\theta), ~\hat{\mathcal{B}}_\infty(\theta) $ and $\hat{\mathcal{C}}_\infty(\theta)$. 
The DFTs $ \hat{\mathcal{A}}_n,~ \hat{\mathcal{B}}_n, ~\hat{\mathcal{C}}_n$ of $ \mathcal{A},~ \mathcal{B}, ~\mathcal{C}$ are now precisely the values of $ \hat{\mathcal{A}}_\infty(\theta), ~\hat{\mathcal{B}}_\infty(\theta) $ and $\hat{\mathcal{C}}_\infty(\theta)$ at $\theta = \frac{2\pi}{L}n$, for all wavenumbers $n \in \Znd$. 

\subsection{ \hn norm evaluation in the spatial frequency domain}
From now on, let us assume that the system~\eqref{eq:generalSS} is input-output stable, so that its \hn norm \eqref{eq:variancedef} exists. This norm can then be calculated as
\begin{equation}
\label{eq:h2normclassic}
\mathbf{V}_N = \mathrm{tr}\left( \int_0^\infty \mathcal{B}^*e^{\mathcal{A}^*t}\mathcal{C}^*\mathcal{C}e^{\mathcal{A}t}\mathcal{B} \mathrm{d}t \right).
\end{equation}
Now, recall that the system~\eqref{eq:generalSS} could be (block) diagonalized by the DFT, where the Fourier symbols $\hat{\mathcal{A}}_n,\hat{\mathcal{B}}_n, \hat{\mathcal{C}}_n$ correspond to the decoupled diagonal elements. Since the DFT is a unitary transformation towards which the \hn norm is invariant, the trace in~\eqref{eq:h2normclassic} can be re-written as:
\begin{equation}
\label{eq:normexpression}
\mathbf{V}_N  = \mathrm{tr} \left( \sum_{n \in \Znd} \int_0^\infty \hat{\mathcal{B}}^*_n e^{\hat{\mathcal{A}}_n^*t} \hat{\mathcal{C}}^*_n \hat{\mathcal{C}}_n e^{\hat{\mathcal{A}}_nt} \hat{\mathcal{B}}_n \mathrm{d}t \right)
\end{equation}
Now, consider the output operator $H$ defined in \eqref{eq:Hdavdef}. It is easy to verify that its Fourier symbol is $\hat{h}_0 = 0$, and $\hat{h}_n = 1$ for $n \neq 0$. This implies that the output matrix $\hat{\mathcal{C}}_0 = 0$ for all systems considered in this paper (i.e., the zero mode is unobservable). Consequently, we can obtain the \hn norm in~\eqref{eq:normexpression} by summing only over $n \in \Znd \backslash \{0\}$.

Furthermore, following the discussion in the previous section, we can regard the Fourier symbols in \eqref{eq:normexpression} as subsamples of $\hat{\mathcal{A}}_\infty(\theta), \hat{\mathcal{B}}_\infty(\theta)$, and $\hat{\mathcal{C}}_\infty(\theta)$. Given this relationship, we can now state the per-site variance $\Vn=\mathbf{V}_N/N$ from~\eqref{eq:pernodevariance} as 
\begin{equation} \label{eq:samplingsum}
\Vn = \frac{1}{N} \sum_{\mathclap{\substack{ \theta = \frac{2\pi}{L}n \\ n \in \Znd \backslash \{0\}} } } \mathrm{tr} \left( \hat{\mathcal{B}}^*_\infty (\theta) \hat{P}(\theta) \hat{\mathcal{B}}_\infty(\theta) \right), \end{equation}
where the individual time integrals are defined as follows:
\begin{definition}
\label{def:gramian}
\begin{equation}
\label{eq:gramiandef}
{\hat{P}(\theta)} := \int_0^\infty e^{\hat{\mathcal{A}}_\infty^*(\theta)t} \hat{\mathcal{C}}^*_\infty (\theta) \hat{\mathcal{C}}_\infty (\theta)e^{\hat{\mathcal{A}}_\infty(\theta) t} \mathrm{d}t.  
\end{equation}
We call $\hat{P}(\theta)$ the \emph{observability Gramian} at $\theta$. 
\end{definition} 
%\vspace{2mm}
The observability Gramian at each $\theta \neq 0$ is obtained by solving the Lyapunov equation
\begin{equation}
\label{eq:lyap}
\hat{\mathcal{A}}^*_\infty(\theta)\hat{P}(\theta)+\hat{P}(\theta)  \hat{\mathcal{A}}_\infty(\theta) = -\hat{\mathcal{C}}^*_\infty(\theta)\hat{\mathcal{C}}_\infty(\theta),
\end{equation}
and is unique and finite provided $\hat{\mathcal{A}}_\infty(\theta)$ is Hurwitz. 

For all problem formulations considered here, $\hat{\mathcal{B}}_\infty(\theta)$ is a vector where one element\footnote{In the vehicular formation case, each ``element'' is a $d\times d$ diagonal matrix, where each of the $d$ diagonal elements is equal by Assumption~\ref{ass:decoupling}.} is 1 and remaining elements are zero. Thus, $\mathrm{tr}( \hat{\mathcal{B}}^*_\infty(\theta) \hat{P}(\theta) \hat{\mathcal{B}}_\infty(\theta) )$ in \eqref{eq:samplingsum} is just one element of the matrix $\hat{P}(\theta)$\footnote{Or the sum of $d$ identical such elements.}. This is a quantity that will be used throughout the paper and we make the following definition:
\begin{definition}[Per-site \hnnorm density]
\label{def:trgramian}
\begin{equation}
\label{eq:trgramian}
{\hat{p}(\theta)} := \mathrm{tr}\left( \hat{\mathcal{B}}^*_\infty(\theta) \hat{P}(\theta) \hat{\mathcal{B}}_\infty(\theta) \right).
\end{equation}
This quantity captures the distribution of the per-site variance~$\Vn$ over the spatial frequency variable~$\theta$ and we therefore refer to it as the \emph{per-site \hnnorm density}.
\end{definition} \vspace{0.5mm}

Now, notice that if the value of $\hat{p}(\theta)$ is bounded for all 
%$\theta \in [-\pi,\pi]^d$, 
$\theta \in \mathcal{R}^d$, 
then $\Vn$ in \eqref{eq:samplingsum} will remain bounded as $N \rightarrow \infty$ and the system in question is to be regarded as \emph{fully coherent}. For the consensus and vehicular formation problems, however, there is typically a single zero eigenvalue at wavenumber $n =0$ that corresponds to the spatial average mode (see Section~\ref{sec:perfmeas}). This makes $\hat{\mathcal{A}}_\infty(0)$ non-Hurwitz, and in turn causes a singularity in $\hat{p}(\theta)$ at $\theta = 0$. Even though the mode at $\theta =0$ itself is unobservable from the system output, the singularity makes the \hnnorm density grow unboundedly for small $\theta$, that is, for small wavenumbers.

For this reason, we use the following appropriate integral to estimate the value of the sum in \eqref{eq:samplingsum}:% asymptotically:
\begin{equation}
\label{eq:integraldef}
{S(\Delta)} := \int_{\Delta\le |\theta_1| \le \pi} \cdots  \int_{\Delta\le |\theta_d| \le \pi}   \hat{p}(\theta)~\mathrm{d}\theta_1 \cdots \mathrm{d}\theta_d ,
\end{equation}
where the argument $\Delta$ indicates the size of a deleted neighborhood around $\theta = 0$.  We recognize the sum in \eqref{eq:samplingsum} as a Riemann sum approximation of the integral \eqref{eq:integraldef}  with volume element $1/N = 1/L^d$. The integral can therefore be used to bound the sum asymptotically. Consider the following lemma:
\begin{lemma}
\label{thm:riemannapprox}
The per-site variance $\Vn$ in \eqref{eq:samplingsum} is upper and lower bounded by the integral \eqref{eq:integraldef} as
\begin{equation}
\label{eq:riemannbounds}
S\left( \frac{4\pi}{L} \right) \le \Vn \le S\left(\frac{2\pi}{L} \right),
\end{equation} for all $L>\bar{L}$, for some fixed $\bar{L}$.
\end{lemma}
\begin{IEEEproof} 
See appendix.
\end{IEEEproof}
The integral and the Riemann sum approximations are illustrated in Figure~\ref{fig:riemann}.

\begin{figure}
\centering
\includegraphics[width = 0.46\textwidth]{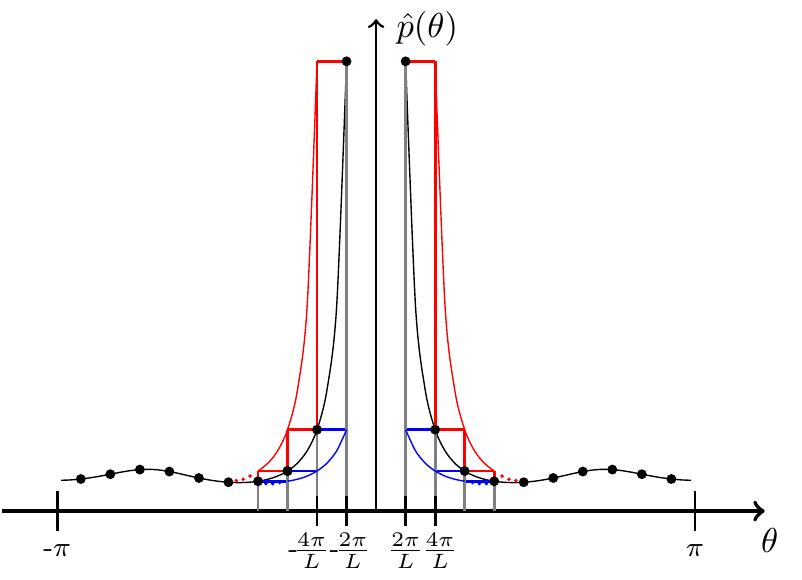}
\caption{Illustration of the upper and lower bounds in Lemma~\ref{thm:riemannapprox}. 
The per-site variance $V_N$ is a sum that can be bounded by the lower (upper) Riemann integrals of the \hnnorm density $\hat{p}(\theta)$ {shifted} by $\frac{2\pi}{L}$  ($\frac{4\pi}{L}$) represented in red (blue). 
{The systems we consider have a zero eigenvalue at $\theta = 0$, causing a singularity in $\hat{p}(\theta)$.  The order of this singularity, i.e., the rate at which $ \hat{p}(\theta)\stackrel{\theta\rightarrow 0}{\longrightarrow} \infty$ determines how fast the integrals, and thereby $\Vn$, grows as $L\rightarrow \infty$. This growth corresponds to the asymptotic performance scaling.  }
%Thus, the rate with which $\hat{p}(\theta) \stackrel{\theta\rightarrow 0}{\longrightarrow} \infty$  determines the asymptotic scaling of $\Vn$.
   }
\label{fig:riemann}
\end{figure}

The performance of the consensus and vehicular formation systems can now be evaluated as follows. First, the system operators are re-defined on $\mathbb{Z}^d$ and (block) diagonalized using the $Z$-transform \eqref{eq:ZtransformPrel}. Second, the Lyapunov equation \eqref{eq:lyap} is solved to determine $\xtr$. Bounds on the variances~$\Vn$ are then found through Lemma~\ref{thm:riemannapprox}. Next, we derive general expressions for the scaling (in $L$) of the integral~\eqref{eq:integraldef}.

\begin{rem}
{It is important to note that the systems we consider remain of \emph{finite} size~$N$ throughout. This is a preqrequisite for the finiteness of the  \hn norm. Only the system operators are re-defined onto the infinite lattice, to facilitate an estimation of the \hn norm by the integrated \hnnorm density according to Lemma~\ref{thm:riemannapprox}. }
\end{rem}
\subsection{Bounds on asymptotic scalings }
\label{sec:boundsgeneral}
We are interested in the scaling of the per-site variance~$\Vn$ in \eqref{eq:samplingsum} with the total number of nodes~$N$ as this number grows large. Using the integral in \eqref{eq:integraldef} and the bounds in Lemma~\ref{thm:riemannapprox}, we can now derive asymptotic scalings of $\Vn$ by exploiting bounds on the per-site \hnnorm density~$\xtr$.
{We begin by a simple example.}

\begin{example}
\label{ex:cons}
Consider the standard consensus algorithm \eqref{eq:standardconsensus} and for simplicity let the dimension $d = 1$. The Lyapunov equation \eqref{eq:lyap} is scalar and solved by
\begin{equation} \label{eq:gramianconsensus1}
\hat{P}(\theta) = \hat{p}(\theta) = \frac{1}{2}\frac{-1}{\hat{f}_\infty(\theta)} \end{equation}
for all $\theta \in [-\pi,\pi] \backslash \{0\}$, where we have used that $\hat{\mathcal{C}}_\infty(\theta) = \hat{h}_\infty (\theta) = 1$ for $\theta \neq 0$. The array $f$ was given in \eqref{eq:standardarray} %as $f_0 = -2\tilde{f}$, $f_1 = f_{-1} = \tilde{f}$ and $f_k = 0$ for $|k|>1$, 
and since $f_k = 0$ for $|k|>1$, we construct the corresponding array $f_\infty$ on $\mathbb{Z}$ by letting $k \rightarrow \infty$. Its $Z$-transform \eqref{eq:ZtransformPrel} is:
\begin{equation}
\label{eq:fourierstandard}
\hat{f}_\infty(\theta) = \tilde{f}(-2+e^{j\theta}+e^{-j\theta}) = -2\tilde{f}(1-\cos \theta).
\end{equation}
Substituting this into \eqref{eq:gramianconsensus1}, the integral in \eqref{eq:integraldef} becomes
$ S(\Delta) = \frac{1}{4\tilde{f}} \int_{\Delta\le |\theta| \le \pi}  \frac{1}{1-\cos \theta} \mathrm{d}\theta.$
%Note that the Gramian behaves as $\frac{1}{\theta^2}$ for small $\theta$, making the integral large for small $\Delta$. 
The lower bound in \eqref{eq:riemannbounds} is
\[S\left(\frac{4\pi}{L}\right) = \frac{-1}{2\tilde{f}}\left[ \cot \frac{\theta}{2} \right]_{\frac{4\pi}{L}}^\pi = \frac{1}{2\tilde{f}}\cot \frac{2\pi	}{L}, \]
and the upper bound has the same form.
A series expansion of the cotangent function reveals that this expression scales as $\frac{1}{\tilde{f}} L = \frac{1}{\tilde{f}} N $ asymptotically. 
%The same holds for the upper bound in \eqref{eq:riemannbounds}. 
This result is represented in case 1a) of Theorem~\ref{thm:mainresult} ($\tilde{f}$ here corresponds to the algorithm parameter~$\beta$).%, and is in line with the result in \cite[Table 1]{Bamieh2012}.
\end{example}
\vspace{1mm}
In general, let us assume that the \hnnorm density is such that 
\begin{equation}
\label{eq:dominantbehavior}
\hat{p}(\theta) \sim \frac{1}{\beta^{r/2}} \frac{1}{(\theta_1^2+ \theta_2^2 + \cdots + \theta_d^2)^{r/2}},
\end{equation}
for some non-negative {$r$}. The number $r$ characterizes the order of the \hnnorm density's singularity at $\theta = 0$. %That is, how fast $\hat{p}(\theta) $ grows as $\theta \rightarrow 0$. 
In the upcoming analysis, we will show that any admissible controller for the systems considered in this paper results in \hnnorm densities that satisfy~\eqref{eq:dominantbehavior} with $r \in \{0,2,4\}$.

We have also introduced the algorithm parameter $\beta$, which reflects the size of the system's feedback gains (c.f. $\tilde{f}$ in Example~\ref{ex:cons}). In particular, let $\beta := \max\{ ||f||_\infty,||g||_\infty \}$. All feedback array elements, which are bounded by assumption, are then proportional to $\beta$. We show in Section~\ref{sec:controleffort} that the parameter $\beta$ is bounded by the system's total control effort. It can therefore be considered a proxy for control effort.

The number $r$ determines the coherence properties for a given system. If $r=0$, the system is fully coherent. Otherwise, the level of coherence depends on the spatial dimension $d$ of the network. We now state the main result of this section:
\begin{lemma}
\label{thm:scaling}
Assume that the per-site \hnnorm density~$\hat{p}(\theta)$ defined in~\eqref{eq:trgramian} satisfies \eqref{eq:dominantbehavior}. 
The {steady state} per-site variance~\eqref{eq:pernodevariance} then scales asymptotically as
\begin{equation}
\label{eq:scalinggeneral}
\Vn \sim \frac{1}{\beta^{r/2}} \begin{cases}
 L^{r-d}  & ~~\mathrm{if}~ d < r \\
\log L  & ~~\mathrm{if}~ d= r \\
1 & ~~\mathrm{if}~ d >r
 \end{cases}
\end{equation}
up to some constant, which is independent of the lattice size~$L$ and the algorithm parameter $\beta$. 
\end{lemma}
\begin{IEEEproof}
First, substitute the approximation \eqref{eq:dominantbehavior} into the integral  $S(\Delta)$ in~\eqref{eq:integraldef} and denote the resulting integral $\tilde{S}(\Delta)$.
{We transform this to hyperspherical coordinates by defining $\rho = (\theta_1^2 \!+\!  \cdots \! +\! \theta_d^2)^{1/2}$ and the $d-1$ coordinates $\phi_1,\ldots,\phi_{d-2} \in [0,\pi]$ and $\phi_{d-1}\in [0,2\pi]$ which are such that $\theta_i = \rho \cos \phi_i \prod_{j = 1}^{i-1}\sin \phi_j $ for $i = 1,\ldots,d-1$ and $\theta_d = \rho\prod_{j = 1}^{d-1}\sin \phi_j $. We obtain:}

\vspace{-4mm}
{
\begin{small}
\begin{multline} 
\tilde{S}(\Delta) = \int_{\Delta\le |\theta_1| \le \pi} \!  \cdots \! \int_{\Delta\le |\theta_d| \le \pi} \! \frac{1}{\beta^{r/2}}\! \frac{1}{(\theta_1^2 \!+\!  \cdots \! +\! \theta_d^2)^{r/2}}  \mathrm{d}\theta_1  \! \cdots  \! \mathrm{d}\theta_d \\
=\!\!\frac{c_d}{\beta^{r/2}}  \int_\Delta^\pi \!\! \int_0^{2\pi} \!\! \int_0^\pi \! \! \cdots \!\! \int_0^\pi \!\! \frac{1}{\rho^r} \rho^{d \!-\!1} \sin^{d\!-\!2 }\!\phi_1 \! \cdots \! \sin \phi_{d-2}  \mathrm{d}\rho\mathrm{d}\phi_{d\!-\!1} \cdots \mathrm{d}\phi_{1}\\ \label{eq:integraleval}
= \! \frac{c_d}{\beta^{r/2}} \mathcal{S}_{d}\int_{\Delta}^\pi \rho^{d-r-1} \mathrm{d}\rho,
 \end{multline}
\end{small}}
\vspace{-3mm}

\noindent {where $\mathcal{S}_{d}$ is the (generalized) surface area of the $d$-dimensional unit sphere and $c_d$ is a bounded scaling factor arising from integrating over a hypersphere rather than a hypercube. }

Now, by Lemma~\ref{thm:riemannapprox} we know that $\Vn$ is bounded as 
\[ \underline{c} \tilde{S}\left( \frac{4\pi}{L} \right)  \le S \left( \frac{4\pi}{L} \right) \le  \Vn \le  S \left( \frac{2\pi}{L} \right) \le \bar{c} \tilde{S} \left( \frac{2\pi}{L} \right),\] 
for all $L \ge \bar{L}$ for some $\bar L$, and with the constants $\underline{c}, \bar{c}$ from {the scaling bounds in}~\eqref{eq:dominantbehavior}.
Substituting for $\Delta$ in \eqref{eq:integraleval} the values $\frac{2\pi}{L}$ and $\frac{4\pi}{L}$ from these upper and lower bounds and defining new constants $\underline{c}', \bar{c}'$, {the solution to the integral} gives that
%\begin{subequations}
\begin{align*}
\Vn & \le \bar{c}'\mathcal{S}_d \frac{1}{\beta^{r/2}} \begin{cases}
\frac{1}{r-d} \pi^{d-r} \left( \left(\frac{L}{2}\right)^{r-d} -1\right) & ~~\mathrm{if}~ d\neq r \\
\log L - \log 2 & ~~\mathrm{if}~ d= r
 \end{cases} \\
 \Vn & \ge \underline{c}'\mathcal{S}_d \frac{1}{\beta^{r/2}} \begin{cases}
\frac{1}{r-d} \pi^{d-r} \left( \left(\frac{L}{4}\right)^{r-d} -1\right) & ~~\mathrm{if}~ d\neq r \\
\log L - \log 4 & ~~\mathrm{if}~ d= r  \end{cases} 
\end{align*}
%\end{subequations}
Noticing that these bounds are identical up to a constant for any given $d$, the result \eqref{eq:scalinggeneral} follows. 
\end{IEEEproof}
\vspace{1mm}

%%%%%%%%%%%%%%%%%%%%%%%%%%%%%%%%%%%%%%%%%%%%%%%%%%%%%%%%%%%%%%%%%%%%%%%%%%%%%%%%%%%
% Stability Section
%%%%%%%%%%%%%%%%%%%%%%%%%%%%%%%%%%%%%%%%%%%%%%%%%%%%%%%%%%%%%%%%%%%%%%%%%%%%%%%%%%%
\section{Admissibility of Dynamic Feedback Laws}
\label{sec:stability}
We now turn to the question of stability of the consensus and vehicular formation systems with dynamic feedback, which is a prerequisite for the \hn performance evaluation laid out in the previous section. In particular, we must require the underlying systems to be BIBO stable for any network size $N$ to allow for the asymptotic performance analysis. 

%Clearly, the systems' stability in the limit of an infinite network is a prerequisite for their asymptotic performance evaluation. 

{With} static feedback, BIBO stability can easily be guaranteed by ensuring that the feedback operators $F$ and $G$ in the systems~\eqref{eq:consensusstatic} and~\eqref{eq:vehiclestatic} have negative Fourier symbols, i.e., $\hat{f}_n, ~\hat{g}_n<0$ for all wavenumbers $n \neq 0$. With dynamic feedback, on the other hand, Routh-Hurwitz stability criteria must typically be evaluated on a case-by-case basis to derive sufficient stability conditions. It turns out, however, that certain feedback configurations will inevitably lead to instability beyond a certain network size. Such feedback laws are therefore \emph{inadmissible} with respect to our analysis. In order to rule them out, this section presents necessary conditions for stability { at any network size}.%the limit of large networks.  %of the the consensus and vehicular formation systems with dynamic feedback. 

\subsection{Conditions for input-output stability}
The stability of a given LTI system on the form~\eqref{eq:generalSS} can be verified by ensuring that its individual Fourier symbols are stable in their own right. 
We begin by re-stating the following Theorem from previous work:
\begin{theorem} \cite[Corollary 1]{Bamieh2002}
\label{thm:stability}
The system \eqref{eq:generalSS} on $\Znd$ is exponentially stable if and only if the matrix~$\hat{\mathcal{A}}_n$ is Hurwitz stable for every $n \in \Znd$. 
\end{theorem}
\begin{IEEEproof}
See \cite[Theorem 1]{Bamieh2002} and note that the group $\Znd$ is compact. 
\end{IEEEproof}
Now, we are evaluating these systems asymptotically, and must therefore require that they remain stable for \emph{any} lattice size $L$, as this number grows. Since the Fourier symbols~$\hat{\mathcal{A}}_n$ can be seen as subsamples of $\hat{\mathcal{A}}_\infty(\theta)$ (see Section~\ref{sec:operatorsinfinite}), the only way to ensure stability for any lattice size $L$ is to make sure that $\hat{\mathcal{A}}_\infty(\theta)$ is stable for every~$\theta$. In our case, though, the mode at $n = 0$ is unobservable from the considered output (see Remark~\ref{rem:BIBO}). BIBO stability is therefore guaranteed if $\hat{\mathcal{A}}_\infty(\theta)$ is stable for every~$\theta$ away from zero:
\begin{corollary}
\label{thm:stabilitycor}
The system \eqref{eq:generalSS} on $\Znd$ with output defined as in~\eqref{eq:Hdavdef} {is} BIBO stable for any network size~$N = L^d$ if and only if the matrix $\hat{\mathcal{A}}_\infty(\theta)$ is Hurwitz stable for all~$\theta \in\mathcal{R}^d\backslash \{0\}$.
\end{corollary}

In order to constrain the upcoming performance analysis to feedback laws that guarantee stability for any network size~$N$ according to Corollary~\ref{thm:stabilitycor}, we make the following definition:
\begin{definition}[Admissible feedback law]
A feedback control law defined on~$\Znd$ is \emph{admissible} if and only if the corresponding closed-loop system is BIBO stable with respect to the output~\eqref{eq:Hdavdef} for any network size $N = L^d$.
\end{definition}
\begin{rem}
\label{rem:BIBOstability}
Note that the considered systems are finite-dimensional for any given lattice size $L$. Their BIBO stability is therefore equivalent to the total variance~$\mathbf{V}_N$ being bounded.%, see also Remark~\ref{rem:BIBO}.
\end{rem}

{Under relative feedback, admissibility of the dynamic feedback laws is not straightforward.
We next present some necessary conditions.}

\subsection{Admissibility conditions under relative feedback}
First, consider the consensus problem with dynamic feedback~\eqref{eq:consensusdynamic} with feedback operators $A,B,F$. Using Corollary~\ref{thm:stabilitycor}, we derive the following theorem:
\begin{theorem}
\label{thm:newstabilitycons}
Consider the consensus system~\eqref{eq:consensusdynamic}. The feedback law is admissible \emph{only if} at least one of the following conditions holds:
%The system \eqref{eq:consensusdynamic} can be input-output stable with respect to the output \eqref{eq:Hdavdef} for any lattice size $N$ \emph{only if} at least one of the following conditions holds:
\begin{enumerate}[a)]
\item The operator $B$ is \emph{symmetric}, %implying that $\hat{b}_\infty(\theta)$ is real,
\item The operator $A$ involves \emph{absolute feedback}, that is, $A$ does not satisfy \eqref{eq:relativedef}.
\end{enumerate}
\end{theorem} 
\begin{IEEEproof}
See appendix. 
\end{IEEEproof}

{In the vehicular formation case with relative velocity feedback, a similar admissibility condition holds: }
\begin{theorem}
\label{thm:stabilityvehicle1} 
Consider the vehicular formation system~\eqref{eq:vehicledynamic}, where the feedback operators $F,G,B,C$ have the relative measurement property~\eqref{eq:relativedef}. The feedback law is admissible \emph{only if} at least one of the following conditions holds:
\begin{enumerate}[a)]
\item The operator $B = 0$, while $A \neq 0$,
\item The operator $A$ involves \emph{absolute feedback}, that is, $A$ does not satisfy \eqref{eq:relativedef}.
\end{enumerate}
\end{theorem}
\begin{IEEEproof}
See appendix. 
\end{IEEEproof}
Theorem~\ref{thm:stabilityvehicle1} implies that integral control based on position measurements cannot be implemented for large networks, unless there is %built-in ``damping'' of the auxiliary state through 
an absolute feedback term in~$A$. 
{ Note, however, that if the purpose of the dynamic feedback law is to eliminate stationary  errors through integral action, including such a term in~$A$ would defeat the purpose. In this case, the auxiliary state $z$ is namely stabilized, and the integral action reduced to zero. }
\begin{rem}
\label{rem:stabilityrem} 
Theorems~\ref{thm:newstabilitycons}--\ref{thm:stabilityvehicle1} imply that a system with a given feedback protocol may be stable for small lattice sizes~$L$, but becomes unstable at some lattice size $L^{\mathrm{crit}}$ unless the criteria are satisfied. As long as the control effort (feedback gains) is bounded, $L^{\mathrm{crit}}$ will always exist and be finite.
\end{rem}

%%%%%%%%%%%%%%%%%%%%%%%%%%%%%%%%%%%%%%%%%%%%%%%%%%%%%%%%%%%%%%%%%%%%%%%%%%%%%%%%%%%
% Evaluation
%%%%%%%%%%%%%%%%%%%%%%%%%%%%%%%%%%%%%%%%%%%%%%%%%%%%%%%%%%%%%%%%%%%%%%%%%%%%%%%%%%%

\section{Performance Scalings with Dynamic Feedback}
\label{sec:evaluation} 
We established in Section~\ref{sec:finiteinfinite} that the asymptotic performance scaling depends on properties of the per-site \hnnorm density. We {now} evaluate the \hn norm densities for admissible feedback laws and derive this paper's main result that was previewed in Theorem~\ref{thm:mainresult}. {In order to establish results for dynamic feedback, we first need to consider the respective problem under static feedback. }

\subsection{Consensus: performance with static feedback}
\label{sec:conseval}
We begin by deriving the performance scaling for the static consensus problem~\eqref{eq:consensusstaticcontrol}. As seen in Example~\ref{ex:cons}, the Lyapunov equation \eqref{eq:lyap} is a scalar equation, which is solved by 
\begin{equation}
\label{eq:gramianconsensus}
\hat{P}(\theta) = \hat{p}(\theta)=  \frac{-1}{2\mathrm{Re}\{\hat{f}_\infty (\theta)\}} .
\end{equation}
%To characterize this \hnnorm density, we use the following {lemma}:
{Now, consider the following lemma:}
\begin{lemma} 
\label{prop:Fscale}
Consider the static consensus system~\eqref{eq:consensusstaticcontrol}. Provided the feedback operator~$F$ is admissible, it holds
%feedback operator~$F$ in the static consensus problemand let Assumptions~\ref{ass:spatialinvariance}--\ref{ass:relative} hold. Provided $F$ is admissible, it holds
\begin{equation}
\label{eq:fscaling}
\mathrm{Re}\{\hat{f}_\infty (\theta)\} \sim -\beta (\theta_1^2 + \ldots + \theta_d^2 ).
\end{equation}
\end{lemma}
\begin{IEEEproof}
By the definition of the $Z$-transform \eqref{eq:ZtransformPrel} it holds
\begin{multline}\label{eq:fouriercos}
\mathrm{Re}\{\hat{f}_\infty (\theta)\}  = \!\! \sum_{k \in \mathbb{Z}^d} \!\! f_k \cos(\theta \! \cdot \! k) = \!\!\sum_{k \in \mathbb{Z}^d}\!\! f_k \left[ 1 \!- \! \left(1-\cos(\theta \! \cdot \! k )\right) \right]   \\
  = -\!\! \sum_{k \in \mathbb{Z}^d}\!\! f_k  \left(1-\cos(\theta \! \cdot \! k )\right),
\end{multline} 
where we have used the relative measurement property~\eqref{eq:relativedef}, which implies $\sum_{k \in \mathbb{Z}^d} f_k = \sum_{k \in \mathbb{Z}^d_N}f_k = 0$. A Taylor series expansion of~\eqref{eq:fouriercos} around $\theta = 0$ is
\begin{align}
 \sum_{k \in \mathbb{Z}^d} \!\! f_k  \! \left( \! 1 \! - \! \cos(\theta \!\cdot\! k )\right) = \!\! \sum_{k \in \mathbb{Z}^d} \!\! f_k \! \left( \! \frac{(\theta \!\cdot\! k )^2}{2} \!-\! \frac{(\theta \!\cdot\! k )^4}{4!} \! +\! \cdots \! \right) \!,
 \label{eq:taylor}
\end{align}
which is upper bounded by its first term, that is $\sum_{k \in \mathbb{Z}^d} f_k \! \left(1 \!-\!\cos(\theta \!\cdot\! k )\right) \! \le \! \sum_{k \in \mathbb{Z}^d} f_k   \frac{(\theta\cdot k)^2}{2} $ for all~$\theta$. Thus,
\begin{multline}
  \sum_{k \in \mathbb{Z}^d} f_k  \left(1-\cos(\theta\cdot k )\right)  \le  \frac{1}{2} \sum_{k \in \mathbb{Z}^d} |f_k| \left( k_1\theta_1  + \cdots +k_d\theta_d \right)^2 \\
\le \frac{1}{2}\sum_{0 \neq k \in \mathbb{Z}^d} ||f||_\infty q^2 \left( |\theta_1| + \cdots +|\theta_d| \right)^2 \\ 
 \le 2^{d-1} q^{d+2} ||f||_\infty (2d+1)(\theta_1^2 + \cdots +\theta_d^2),  \label{eq:lowerboundgramian}
\end{multline}
where the second inequality follows from the locality assumption~\eqref{eq:localitycondition} and the third from straightforward algebra.

Next, the Taylor expansion~\eqref{eq:taylor} reveals that $\mathrm{Re}\{\hat{f}_\infty (\theta)\}$ goes to zero at a quadratic rate. We can therefore always find a fixed, nonnegative $\underline{c}$ so that 
\begin{equation}
-\mathrm{Re}\{\hat{f}_\infty (\theta)\} =\sum_{k \in \mathbb{Z}^d} f_k  \left(1-\cos(\theta\cdot k )\right) \ge \underline{c}(\theta_1^2 + \cdots +\theta_d^2) \label{eq:lowerboundfourier}
\end{equation}
in some interval near zero; $\theta \in [-\Delta,\Delta]^d$ for a small~$\Delta$. Note that no lower-degree polynomial in~$\theta$ (apart from the zero polynomial) could serve as a lower bound in~\eqref{eq:lowerboundfourier}.
Furthermore, given that the feedback law is admissible, it must hold $-\mathrm{Re}\{\hat{f}_\infty (\theta)\}>\epsilon$ for all $\theta \in \mathcal{R}^d \backslash  [-\Delta,\Delta]^d  =: \mathcal{R}_\Delta^d$ with any fixed $\Delta$. We can therefore always adjust $\underline{c}$ so that~\eqref{eq:lowerboundfourier} holds for the entire region~$\mathcal{R}^d$.

Defining the algorithm parameter $\beta=||f||_{\infty}$, and noticing that remaining parameters of~\eqref{eq:lowerboundgramian} and \eqref{eq:lowerboundfourier} are independent of~$\theta$ and~$L$, the result~\eqref{eq:fscaling} follows. 
\end{IEEEproof}

Inserting the scaling from Lemma~\ref{prop:Fscale} into the \hnnorm density~\eqref{eq:gramianconsensus} shows that
\begin{equation}
\label{eq:scalingconsensus2}
\hat{p}(\theta) =  \frac{-1}{2\mathrm{Re}\{\hat{f}_\infty (\theta)\}} \sim \frac{1}{\beta} \frac{1}{(\theta_1^2 + \cdots +\theta_d^2) },
\end{equation}
that is, the \hnnorm density for the static consensus system~\eqref{eq:consensusstaticcontrol} satisfies \eqref{eq:dominantbehavior} with $r=2$. The per-site variance thus scales according to Lemma~\ref{thm:scaling} with $r =2$. 

\subsection{Consensus: performance with dynamic feedback}
{Before turning to the case of dynamic feedback,} note that the performance of the consensus system with static feedback is independent of any imaginary part of the Fourier symbol~$\hat{f}_\infty(\theta)$. It is therefore independent of whether the feedback operator $F$ is symmetric or not. In the upcoming evaluation of dynamic feedback, we will therefore limit the analysis to $F$ being symmetric:
\begin{ass}
\label{ass:F}
The operator $F$ in the dynamic consensus protocol \eqref{eq:consensusdynamic} is symmetric, that is, it satisfies the properties listed in Assumption~\ref{ass:symmetry}. {It follows that $\hat{f}_\infty(\theta) = \mathrm{Re}\{\hat{f}_\infty(\theta)\}$.}
\end{ass}
\begin{rem}
Assumption~\ref{ass:F} is made to simplify the exposition by limiting the number of possible feedback configurations that must be considered. It is our belief, based on computer-aided evaluation, that the main result would hold also without this assumption. 
\end{rem}

%We now turn to the case of dynamic feedback \eqref{eq:consensusdynamic}, and 
{Let us now} assume that the choice of operators~$A,B,F$ is admissible. % that is, that the stability conditions from Section~\ref{sec:stability} are fulfilled. 
The solution to the Lyapunov equation \eqref{eq:lyap} then gives that
\begin{equation}
\label{eq:gramianfirst2}
\hat{p}(\theta)  = \frac{ -1  }{ 2\hat{f}_\infty(\theta) + 2\varphi^c(\theta)  },
\end{equation}
where $\varphi^c(\theta)$ is a function of the Fourier symbols of $A,B$ and $F$. This \hnnorm density would scale different from \eqref{eq:scalingconsensus2} if the function $\varphi^c(\theta)$ were non-zero and \emph{scaled differently} in~$\theta$ than $\hat{f}_\infty(\theta)$, for which we established { Lemma~\ref{prop:Fscale}}. This is, however, not the case for any admissible configuration of the feedback operators $A$ and $B$. Consider the following lemma:
\begin{lemma}
\label{lem:qscalingcons}
For any admissible choice of the operators $A,B,F$ in \eqref{eq:consensusdynamic} satisfying Assumptions~\ref{ass:spatialinvariance}--\ref{ass:relative}, \ref{ass:F}, the function~$\varphi^c(\theta)$ in \eqref{eq:gramianfirst2} is such that
\begin{equation}
\label{eq:qscaling}
\hat{f}_\infty(\theta) + \varphi^c(\theta)  \sim - \beta (\theta_1^2 + \ldots + \theta_d^2 ).
\end{equation}
%If $A = 0$, it holds that $\varphi^c(\hat{a}_\infty(\theta), \hat{f}_\infty(\theta), \hat{b}_\infty(\theta))  \equiv 0$.
Therefore, the \hnnorm density~$\hat{p}(\theta)$ in \eqref{eq:gramianfirst2} will satisfy \eqref{eq:dominantbehavior} with $r=2$ for any design of the dynamic feedback law.
\end{lemma}
\begin{IEEEproof}
See appendix.
\end{IEEEproof}
The asymptotic performance scaling will thus be unchanged compared to static feedback. Rewriting the asymptotic scalings from Lemma~\ref{thm:scaling} in terms of total network size $N = L^d$ gives the result in Theorem~\ref{thm:mainresult}.

\subsection{Vehicular formations: performance with static feedback}
\label{sec:perfevalvehichle}
{Consider the vehicular formation problem under static feedback~\eqref{eq:vehiclestatic} }. The solution to the Lyapunov equation \eqref{eq:lyap} gives the \hnnorm density 
\begin{equation}
\label{eq:gramianvehicle}
\hat{p}(\theta) = \frac{d}{2\hat{f}_\infty(\theta) \hat{g}_\infty(\theta)}.
\end{equation}
The following {lemma is used} to bound this \hnnorm density:% for the cases with relative and absolute {velocity} feedback:
\begin{lemma}
\label{prop:fgvehicle}
{Consider the feedback operators $F$ and $G$ in the static vehicular formation problem~\eqref{eq:vehiclestatic}, and assume they are admissible. It holds $\hat{f}_\infty(\theta) \sim -\beta (\theta_1^2 + \ldots + \theta_d^2 )$. If $G$ has the relative measurement property~\eqref{eq:relativedef}, then also  $\hat{g}_\infty(\theta) \sim -\beta (\theta_1^2 + \ldots + \theta_d^2 )$. Otherwise, $\hat{g}_\infty(\theta)\sim \hat{g}_0$, for a given constant~$\hat{g}_0$.}
\end{lemma}
\begin{IEEEproof} By Assumption~\ref{ass:symmetry}, $\hat{f}_\infty(\theta),\hat{g}_\infty(\theta)$ are real valued. If they satisfy the relative measurement property~\eqref{eq:relativedef}, they therefore have the same properties as $\mathrm{Re}\{\hat{f}_\infty (\theta)\}$ from the consensus case, and scale as in~\eqref{eq:fscaling}. {If $G$ has absolute feedback}, it follows from~\eqref{eq:fouriercos} that $\hat{g}_\infty(\theta) = \hat{g}_0  -\sum_{k \in \mathbb{Z}^d}  g_k  \left(1-\cos(\theta \! \cdot \! k )\right)$, where $\hat{g}_0 = \sum_{k \in \mathbb{Z}^d}g_k <0$. Due to the locality assumption~\ref{ass:locality}, this {number} is uniformly bounded for all $\theta \in \mathcal{R}^d$, see \eqref{eq:lowerboundgramian}. We can {thus} write $\hat{g}_\infty(\theta) \sim \hat{g}_0$.
\end{IEEEproof}
In the case of {only \emph{relative} feedback}, Lemma~\ref{prop:fgvehicle} bounds the \hnnorm density {from~\eqref{eq:gramianvehicle} }as 
\begin{equation} \label{eq:vehicledefaultscaling}\hat{p}(\theta) = \frac{d}{2\hat{f}_\infty(\theta) \hat{g}_\infty(\theta)} \sim \frac{1}{\beta^2(\theta_1^2 + \ldots + \theta_d^2 )^2 }.\end{equation}
The per-site variance thus scales as in Lemma~\ref{thm:scaling} with~$r=4$. 

With {\textit{absolute velocity} feedback} we instead get that 
\[\hat{p}(\theta) \sim \frac{1}{\beta \left(\theta_1^2 + \ldots + \theta_d^2 \right)}.\] 
{In this case, the per-site variance thus} scales as in Lemma~\ref{thm:scaling} with~$r=2$. 
%{In this case} performance therefore scales as in Lemma~\ref{thm:scaling} with~$r=2$. 

{We can also note that relaxing Assumption~\ref{ass:relabs} and allowing absolute feedback from both position and velocity would let }$\hat{f}_\infty(\theta) \sim \hat{f}_0$ and $~\hat{g}_\infty(\theta)\sim \hat{g}_0$, making the \hnnorm density~\eqref{eq:gramianvehicle} uniformly bounded in $\theta$. That is, $r = 0$ in Lemma~\ref{thm:scaling} and the system would be fully coherent. 

The results for the static case outlined above, which are in line with those in \cite[Table 1]{Bamieh2012}, are summarized in Theorem~\ref{thm:mainresult}. 

\subsection{Vehicular formations: performance with dynamic feedback}
Now, consider the vehicular formation system with dynamic feedback on the form \eqref{eq:vehicledynamic}. Provided the feedback configuration is admissible, %that is, that the stability conditions from Section~\ref{sec:stability} are fulfilled, 
the Lyapunov equation \eqref{eq:lyap} gives 
\begin{equation}
\label{eq:gramianrelative}
\hat{p}(\theta)  = \frac{ d  }{ 2\hat{f}_\infty(\theta) \hat{g}_\infty(\theta)+ 2\varphi^v(\theta )  },
\end{equation}
where $\varphi^v(\theta )$ is a function of the Fourier symbols of the operators $A,B,C,F$ and $G$. 
{ We now analyze \eqref{eq:gramianrelative} for the case with both relative and absolute velocity feedback.}

\subsubsection{Relative feedback}
In order for the \hnnorm density in~\eqref{eq:gramianrelative} to scale differently from the static case~\eqref{eq:gramianvehicle}, the function $\varphi^v(\theta)$ would need to scale differently in $\theta$ from the product $\hat{f}_\infty(\theta) \hat{g}_\infty(\theta) $, whose scaling was established in~\eqref{eq:vehicledefaultscaling}. This is, however, not possible with only relative feedback: 
\begin{lemma}
\label{thm:phifunction1}
For any admissible choice of the operators $A,B,C,F,G$ in \eqref{eq:vehicledynamic}
with only relative feedback in $B,C,F,G$, the function~$\varphi^v(\theta)$ in \eqref{eq:gramianrelative} is such that
%the denominator in \eqref{eq:gramianrelative} satisfies
\begin{equation}
\label{eq:phiscalingvehicle}
\hat{f}_\infty(\theta) \hat{g}_\infty(\theta)+ \varphi^v(\theta ) \sim \beta^2 (\theta_1^2 + \ldots + \theta_d^2 )^2,
\end{equation}
Therefore, the \hnnorm density $\hat{p}(\theta)$ in~\eqref{eq:gramianrelative} will satisfy \eqref{eq:dominantbehavior} with $r=4$ for any design of the dynamic feedback.
\end{lemma}
\begin{IEEEproof}
See appendix. 
\end{IEEEproof}
 We conclude that in the case of {only} relative feedback, dynamic feedback on the form \eqref{eq:vehicledynamic} \emph{cannot} improve the asymptotic performance scaling compared to static feedback. 
\begin{rem}
Certain choices of $A,B,C,F,G$ in \eqref{eq:vehicledynamic} may appear as though one can achieve $\varphi^v(\theta) \sim -\beta (\theta_1^2 + \ldots + \theta_d^2 ) $, and thereby improve performance.
%make performance scale as in \eqref{eq:scalinggeneral} with $r=2$. 
For example, if $A = 0$, it holds $\varphi^v(\theta ) = \hat{b}_\infty(\theta) + \hat{c}_\infty(\theta)\hat{g}_\infty(\theta)$ and one may wish to set $B$ as the standard consensus operator~\eqref{eq:standardconsensus}. However, by Theorem~\ref{thm:stabilityvehicle1}, such a choice is inadmissible.% as it always de-stabilizes the system at some lattice size $L^{\mathrm{crit}}$. %The performance metric \eqref{eq:pernodevariance} would therefore not be defined asymptotically.
\end{rem} %\vspace{1.5mm}

\subsubsection{Absolute velocity feedback}
\label{sec:dapianalysis}
In this case, we first consider the distributed-averaging proportional-integral (DAPI) controller \eqref{eq:DAPIexample} for the 1-dimensional vehicular platoon. The solution to the Lyapunov equation yields 
\begin{equation}
\label{eq:gramiandapi}
\hat{p}_{\mathrm{DAPI}}(\theta)=  \frac{1}{ 2\hat{f}\hat{g} -2 \frac{\hat{c}\hat{f}(\hat{a} + \hat{g}) }{ \hat{a}^2 + \hat{g}\hat{a} - \hat{f} } }, 
\end{equation}
where we have left out the $\infty-$subscript and the arguments of the individual Fourier symbols for notational compactness. 

Now, $A$ and $F$ in DAPI are standard consensus operators whose Fourier symbols look like~\eqref{eq:fourierstandard}, while $G = -g_oI$ and $C = -c_oI$, which gives $\hat{g}_\infty(\theta) = -g_o,~\hat{c}_\infty(\theta) = -c_o$.
%By Assumption~\ref{ass:symmetry}, $A$ and $F$ in the DAPI algorithm are standard consensus operators with Fourier symbols as in~\eqref{eq:fourierstandard}, while $G = -g_oI$ and $C = -c_oI$ (with $c_o,g_o>0$ to ensure stability). 
Inserting into \eqref{eq:gramiandapi} yields (after some simplifications):
\[ 
\hat{p}_{\mathrm{DAPI}}(\theta) = \frac{1}{4g_of_+(1 \! - \!\cos\theta) \! + \!2 \frac{c_og_of_+ + 2c_of_+a_+(1-\cos\theta) }{ f_+ + a_+g_o +  2a_+^2(1-\cos\theta) } }, 
\]
which recognize as being uniformly bounded in $\theta \in \mathcal{R}^d$. This implies that already the 1-dimensional vehicular platoon with DAPI control is \emph{fully coherent}.  This is in contrast to the static control law, which yields the performance scaling in~\eqref{eq:resultscalingvehbetter}, and therefore requires 3 spatial dimensions to be fully coherent. %This result 

%{With absolute velocity feedback, other designs on the form~\eqref{eq:vehicledynamic}}
If absolute velocity measurements are available, several designs of the dynamic feedback in \eqref{eq:vehicledynamic} 
can be shown to give the same result as the DAPI controller. In particular, $G$ and $C$ can also include relative feedback %feedback from velocity errors with respect to neighboring vehicles 
and $B$ can be non-zero.

The asymptotic performance scalings for the vehicular formation problem with dynamic feedback are summarized in Theorem~\ref{thm:mainresult}, where they have been re-written in terms of total network size $N = L^d$.

\subsection{Control effort bounds}
\label{sec:controleffort}
In the above derivations, we introduced the algorithm parameter $\beta = \max\{||f||_\infty,||g||_\infty\}$. This parameter affects the performance scaling, as evident from our main result in Theorem~\ref{thm:mainresult}. In particular, if $\beta$ were allowed to increase unboundedly, %and at a rate up to~$N^{3/2}$, 
full coherence could be achieved in any spatial dimension. This is not feasible in any realistic control problem, where the amount of control effort is bounded. We now show that the size of the feedback array elements and therefore $\beta$ are bounded by the total control effort at each network site, which we quantify through:
\begin{equation}
\label{eq:controlvariance}
\mathbb{E}\{u_k^*u_k\},
\end{equation}
that is, the steady state variance of the control signal at each network site. 
In \cite[Lemma 5.1]{Bamieh2012}, such bounds are presented for the case of static feedback. Here, we present bounds for the dynamic feedback case, but limit the analysis to the consensus algorithm with symmetric feedback for the sake of brevity:
\begin{lemma}
\label{lem:controleffort}
Consider the consensus problem with dynamic feedback \eqref{eq:consensusdynamic}, where the feedback operators $A,B,F$ satisfy Assumptions~\ref{ass:locality} and \ref{ass:symmetry}. The following bounds hold: 
\begin{subequations}
\begin{align} \label{eq:finfbound}
\mathbb{E}\{u_k^*u_k\} &\ge \frac{1}{2}||f||_\infty \\ \label{eq:abinfbound}
\mathbb{E}\{u_k^*u_k\} & \ge \sqrt{ \left( \frac{ ||a||_\infty}{4}\right)^2 +\frac{||b||_\infty }{4(2q)^d} } - \frac{||a||_\infty}{4} 
\end{align}
\end{subequations}
\end{lemma}
\begin{IEEEproof}
See appendix.
\end{IEEEproof}
Note that the constants in the bounds are independent of network size. Since we have set $\beta = ||f||_\infty$ and $||a||_\infty, ||b||_\infty \sim \beta$, we can conclude that the asymptotic scalings for the consensus problem in Theorem~\ref{thm:mainresult} will apply to any algorithm with control effort constraints. 
 
%%%%%%%%%%%%%%%%%%%%%%%%%%%%%%%%%%%%%%%%%%%%%%%%%%%%%%%%%%%%%%%%%%%%%%%%%%%%%%%%%%%
% Examples
%%%%%%%%%%%%%%%%%%%%%%%%%%%%%%%%%%%%%%%%%%%%%%%%%%%%%%%%%%%%%%%%%%%%%%%%%%%%%%%%%%%

\section{Implications and Numerical Example}
\label{sec:examples} 
% {Thus far, we have presented limitations to the performance of distributed static and dynamic feedback in toric lattices. This was done in terms of the scaling of global \hn performance, with respect to an output that captures nodal fluctuations with respect to  the network average. 
\begin{figure*}
  \centering
  \subfloat[][{Static feedback, $N=20$}]{
  \includegraphics[width = 0.30\textwidth]{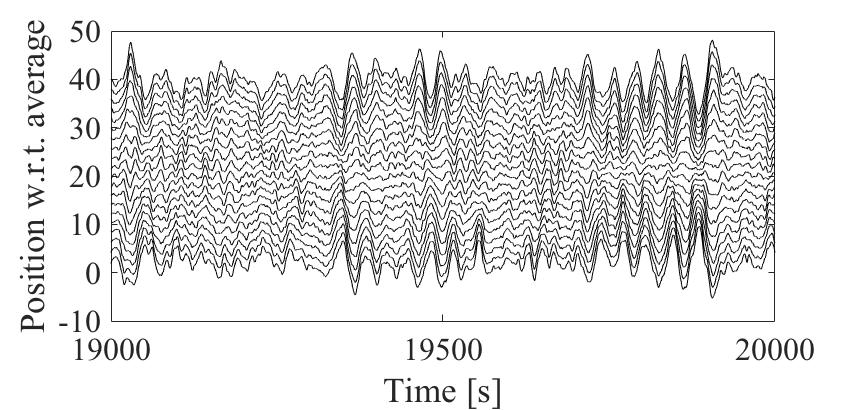}
  \label{fig:stat20}  }
  \subfloat[][{Static feedback, $N=100$}] {
  \includegraphics[width = 0.30\textwidth]{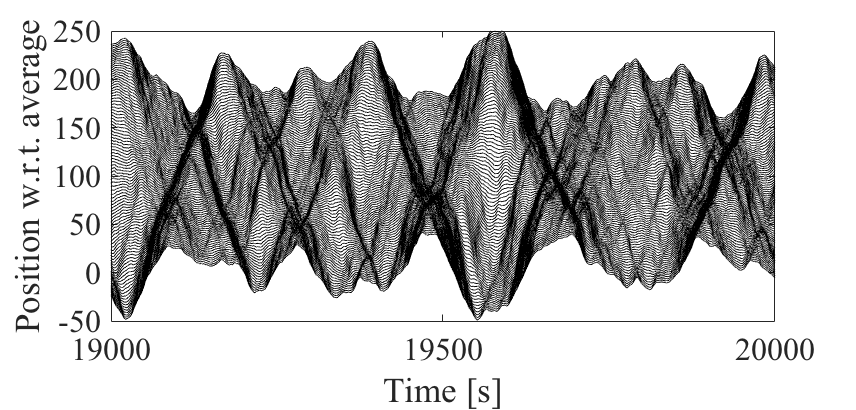}
  \label{fig:stat100}}
   \subfloat[][{Static feedback, $N=200$}] {
  \includegraphics[width = 0.30\textwidth]{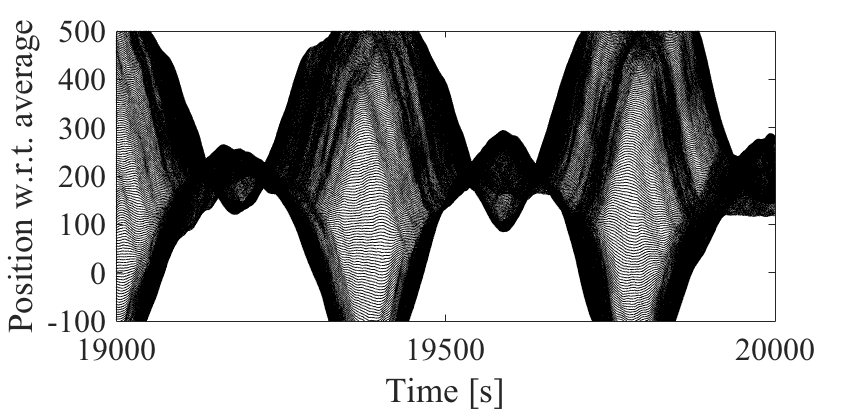}
  \label{fig:stat200}} \\
    \subfloat[][{Dynamic feedback, $N=20$}]{
  \includegraphics[width = 0.30\textwidth]{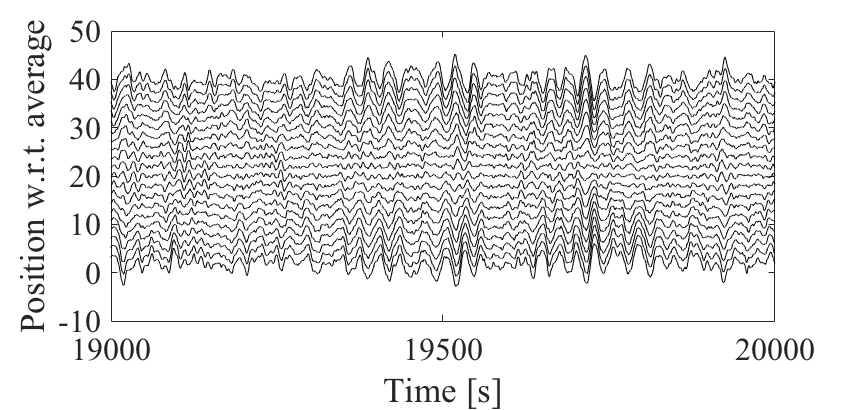}
  \label{fig:dyn20}  }
  \subfloat[][{Dynamic feedback, $N=100$}] {
  \includegraphics[width = 0.30\textwidth]{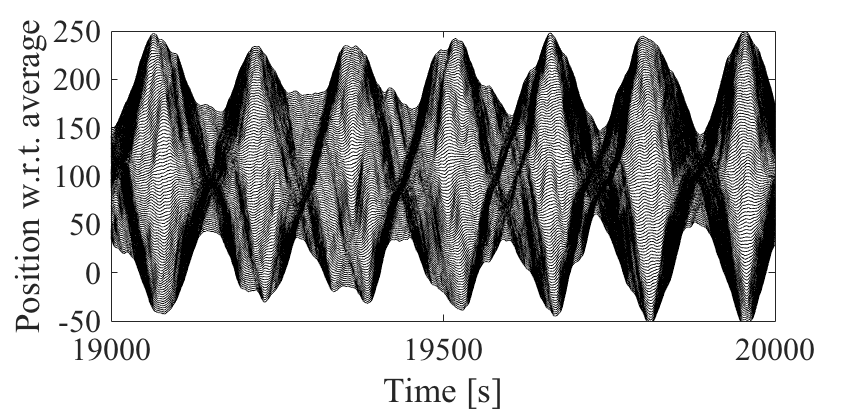}
  \label{fig:dyn100}}
   \subfloat[][{Dynamic feedback, $N=200$}] {
  \includegraphics[width = 0.30\textwidth]{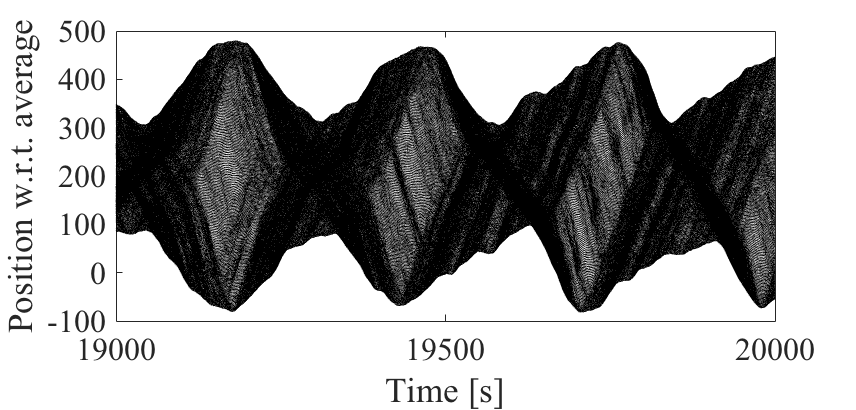}
    \label{fig:dyn200}} \\
      \caption{
{Simulation of an $N$-vehicle platoon with static feedback~\eqref{eq:vehiclestatic} and dynamic feedback~\eqref{eq:vehicledynamic} from relative measurements. At each time step of 0.1 s the independent inputs $w_k$ are sampled from a Gaussian distribution. We display the time trajectories of all vehicles' positions, with the average motion of the platoon subtracted and a reference spacing $\Delta_x = 2$ units inserted between vehicles.  Under perfect control, the trajectories would be $N$ straight horizontal lines separated by $\Delta_x$. Note that the times displayed are $19000~\mathrm{s} \le t \le 20000~\mathrm{s}$ (approx. steady state), and that the scales on the vertical axes are proportional to~$N$. The platoon exhibits an accordion-like motion for large~$N$ with both static and dynamic feedback, showcasing the lack of coherence predicted by Theorem~\ref{thm:mainresult}.}}
	%\vspace*{-1.5em}
\label{fig:simulation}
\end{figure*}

\begin{figure}
\centering
\includegraphics[scale=0.95]{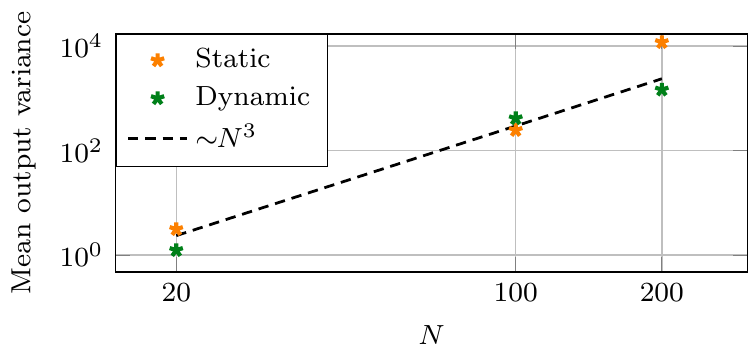}
\caption{{Mean variance (over the platoon) of the performance output~\eqref{eq:performancedef} for the system trajectories displayed in Figure~\ref{fig:simulation}. The data points agree with the $V_N \sim N^3$ scaling predicted by Theorem~\ref{thm:mainresult}.}} %Note the logarithmic axes.}}
\label{fig:scaling}
\end{figure} 

{ The performance limitations discussed in this paper are in terms of scalings of global \hn performance, with respect to an output defined through nodal state fluctuations.
%We have presented limitations of distributed static and dynamic feedback control in terms of scalings of global \hn performance, with respect to an output defined as nodal fluctuations from the network average.
 We argued that a better scaling implies that the network remains more coherent, or rigid, when subjected to a process noise disturbance. Fig.~\ref{fig:simulation} shows simulations of strings of vehicles (i.e. platoons) with both static and dynamic feedback from relative measurements. As the platoons grow, they exhibit an increasing lack of coherence. This is manifested through slow and large-scale fluctuations of the platoon length, clearly indicating that the platoon does not move like a rigid body. While the shape and size of these fluctuations are different with dynamic feedback compared to static, the relative performance deterioration is similar, as predicted by Case 2a in Theorem~\ref{thm:mainresult}. This can also be seen from the corresponding variances displayed in Fig.~\ref{fig:scaling}.  }

{The origin of these slowly varying mode shapes in vehicle platoons was discussed in~\cite{Bamieh2012} and more recently in~\cite{Pates2017}. In this paper, our introduced concept of per-site \hnnorm density provides additional insights. The \hnnorm density is largest near $\theta = 0$, revealing that the low spatial frequency modes are most energetic (see Fig.~\ref{fig:riemann}). As these correspond to the smallest system eigenvalues, they are also temporally slow. %, and when excited not easily attenuated. 
This results in slowly varying modes that have long spatial wavelengths and therefore span across the entire platoon.   }

{The derivations in this paper are made for spatially invariant systems, that is, lattices with periodic boundary conditions. The simulation here, however, is done for a string of vehicles where the first is not connected to the last.  For large platoons, the boundary condition has a limited effect on the interior of the network. The corresponding simulation for a ring of vehicles can indeed be verified to have a very similar appearance to Fig.~\ref{fig:simulation}.}

{While the relation~\eqref{eq:pernodevariance} does not hold if the assumption of spatial invariance is relaxed, the quantity $\Vn$ can be evaluated as the mean variance over the network. This is also what is displayed in Fig.~\ref{fig:scaling}. Through graph embedding (that is, noting that the string can be embedded in a ring graph)~\cite{Barooah2007} it is straightforward to show that the mean variance for the string will be \emph{at least as large} as for the ring graph case. It is therefore subject to the same limitations. Embedding  arguments can also be made in higher spatial dimensions, in particular to relate the performance of 2D lattices to networks described by planar graphs. }  %For general networks, the notions of network scaling, locality, and spatial dimension may however be subtle (see Section~\ref{sec:generalnetwork}.)  } 

{The simulation in Figs.~\ref{fig:simulation}--\ref{fig:scaling} also demonstrates why it is the \emph{scaling} of the per-site variance, rather than its actual value for a given $N$, that is meaningful for describing fundamental limitations. Even though a particular controller can achieve lower per-site variance in a given finite-size network (here, for example, the dynamic one at $N = 200$), the fact that it scales with network size implies that performance inevitably degrades as the network grows. This holds regardless of scaling coefficients. The result of Theorem~\ref{thm:mainresult} thus implies that neither static nor dynamic feedback from only relative measurements is \emph{scalable} to large networks. They are therefore both fundamentally limited.    }

%%%%%%%%%%%%%%%%%%%%%%%%%%%%%%%%%%%%%%%%%%%%%%%%%%%%%%%%%%%%%%%%%%%%%%%%%%%%%%%%%%%
% Discussion
%%%%%%%%%%%%%%%%%%%%%%%%%%%%%%%%%%%%%%%%%%%%%%%%%%%%%%%%%%%%%%%%%%%%%%%%%%%%%%%%%%%

\section{Discussion}

\label{sec:discussion}

\subsection{{Non-regular networks}}

{The results in this paper are derived for systems defined on toric lattice networks, under certain restrictive assumptions. Accepting a generalization of the coherence metric~$\Vn$ in~\eqref{eq:pernodevariance}, the assumptions of symmetry, uniformity in gains, and isotropy can be relaxed at the cost of analytic tractability, or by giving looser bounds on performance. Graph embeddings can, as already discussed, also be used to bound performance of more general networks through the lattices in which they can be embedded. {The principle for this argument is that the removal of any network connection can only decrease the graph Laplacian eigenvalues (corresponding to the Fourier symbols in this study) and therefore increases $V_N$. Any subgraph of a lattice (i.e. any embedded graph) thus has greater $V_N$ than the lattice. See e.g. \cite{Barooah2007,Flamme2018} for details. } We note that our theory allows for~$q$ neighbor connections in each lattice direction, making such embedding arguments less restrictive than they may seem.} 

{Other concepts that are important for this paper's results, such as locality, spatial dimension and a consistent notion of growing the network, are not straightforwardly generalized. For families of graphs where the behavior of the graph Laplacian eigenvalues (i.e., the Laplacian eigenvalue density) is known, the ideas in Section~\ref{sec:finiteinfinite} (e.g. the \hnnorm density) can be applied. The notion of spatial dimension can then likely be generalized to one of \emph{spectral} dimension. However, such considerations would only apply to graph families that can be scaled consistently, preserving properties like locality and dimensionality. Relevant contributions on performance limitations in other types of graphs have been made in \cite{Patterson2014, SiamiMotee2015,Grunberg2016, Pates2017}. A proper generalization of the topological properties that cause the limitations described in this paper, however, remains an open research question. }

\subsection{Performance improvement with distributed integral control}
\label{sec:limitations}

We established that dynamic feedback such as the DAPI algorithm~\eqref{eq:DAPIexample} can yield a fully coherent vehicular formation in any spatial dimension, provided that it has access to absolute measurements of velocities with respect to a global reference frame. This situation is reasonable in actual vehicular platoons, where one can assume that each vehicle's speedometer can provide absolute velocity measurements, while absolute position data, which would have to rely on, for example, GPS is less readily available.

An intuitive explanation to this result, which was also established in~\cite{Tegling2017}, is that the dynamic feedback protocol serves as a distributed integral controller, which integrates absolute measurements of velocities in time to yield a substitute for absolute position data. With absolute data from both position and velocity, formations are known to be fully coherent~\cite{Bamieh2012}. However, as such a strategy is essentially so-called ``dead reckoning'', it can be sensitive to noisy measurements. 

One issue arises when different controllers' memory states~$z_k$ diverge due to slight measurement errors. This issue appears in completely decentralized integral control {and leads to instability}, but can be solved through distributed averaging of the memory states between controllers, see e.g. \cite{Andreasson2014}. In our case, distributed averaging is achieved by {choosing} $A$ in \eqref{eq:vehicledynamic} to be a consensus-type operator, as in the DAPI example~\eqref{eq:DAPIexample}.  

{A second issue is how noise and bias in the} measurements affect  performance. Results on this topic have been reported in~\cite{Tegling2018ARX}, and reveal that the performance improvement {achieved} through DAPI control is highly sensitive to the {design} of the distributed averaging operator~$A$.

\subsection{Higher order dynamic feedback controllers}
The dynamic feedback controllers considered in this paper all contain a single local memory state $z$, as illustrated in Fig.~\ref{fig:structure}. They thus describe a class of distributed proportional-integral (PI) controllers with respect to the system's states. While we show that this type of controller cannot improve performance scalings compared to static, memoryless controllers as long as they are limited to relative state feedback, it is an open question whether a higher-order controller, with an arbitrary number of local states, can. 

Even with a higher number of controller states, however, the limitation to relative state feedback implies that the marginally stable mode at the origin remains, so the Fourier symbol $\hat{\mathcal{A}}(\theta)$ is singular at $\theta = 0$. As a consequence, the \hnnorm density will scale badly near $\theta = 0$ (consider the Lyapunov equation~\eqref{eq:lyap} and note that the right hand side is identity). We therefore conjecture that the unfavorable scaling of performance in low spatial dimensions remains as long as the number of local states is finite.

%%%%%%%%%%%%%%%%%%%%%%%%%%%%%%%%%%%%%%%%%%%%%%%%%%%%%%%%%%%%%%%%%%%%%%%%%%%%%%%%%%%
% Appendix
%%%%%%%%%%%%%%%%%%%%%%%%%%%%%%%%%%%%%%%%%%%%%%%%%%%%%%%%%%%%%%%%%%%%%%%%%%%%%%%%%%%
\appendix

\subsection{Scalings of sums and products}
\label{sec:appScaling}
{Many of the proofs in this appendix are based on the behaviors, or scalings, of functions of Fourier symbols in $\theta$. Here we make some preliminary remarks on such scalings.  }

Recall that the notation $u(\theta) \sim v(\theta)$ implies $\underline{c}v(\theta) \le u(\theta) \le \bar{c}v(\theta)$ for all $\theta \in \mathcal{R}^d = [-\pi,\pi]^d$, where $\underline{c},\bar{c}$ are fixed, positive constants. For example, we write $\hat{f}_\infty(\theta) \sim \beta(\theta_1^2 +\ldots + \theta_d^2)$, or $\hat{f}_\infty(\theta) \sim \beta \theta^2$ for short.

 For products and sums of such functions, it holds $u'(\theta) = u_1(\theta)u_2(\theta) + u_3(\theta) \sim v_1(\theta)v_2(\theta) + v_3(\theta)$, implying that the bounds are $\underline{c}_1 \underline{c}_2 v_1(\theta)v_2(\theta)  + \underline{c}v_3 \le u'(\theta) \le \bar{c}_2v_2(\theta)\bar{c}_2v_2(\theta) + \bar{c}v_3(\theta)$. For a quotient: $u'(\theta) =u_1(\theta) / u_2(\theta) \sim v_1(\theta)/v_2(\theta)$ implies $(\underline{c}_1/\bar{c}_2) v_1(\theta)/v_2(\theta) \le u'(\theta) \le (\bar{c}_1/\underline{c}_2) v_1(\theta)/v_2(\theta)  $.

Therefore, the scalings of functions of Fourier symbols can be determined simply by inserting the individual Fourier symbols' scalings. %to determine scalings of functions of Fourier symbols, one can simply insert the individual Fourier symbols' scalings. 
For example, if $\hat{f}_\infty(\theta) \sim \beta \theta^2, \hat{g}_\infty(\theta) \sim \beta \theta^2$, then $\hat{f}_\infty(\theta) \hat{g}_\infty(\theta) \sim \beta^2\theta^4$ and $\hat{f}_\infty(\theta)/ \hat{g}_\infty(\theta) \sim 1$. This is used throughout to determine scalings of \hn norm densities.

\subsection{Maclaurin expansions of $Z$-transforms}
The Maclaurin series expansions of $Z$-transforms will be used to derive admissibility conditions in Theorems~\ref{thm:newstabilitycons}--\ref{thm:stabilityvehicle1}. Consider an operator $A$, and its $Z$-transform $\hat{a}_\infty(\theta)$ given in~\eqref{eq:ZtransformPrel}. The Maclaurin expansion of $\hat{a}_\infty(\theta)$ in the coordinate direction $\theta = (\theta_1,0,\ldots,0)$ is 
\begin{equation}
\label{eq:maclaurin}
\hat{a}_\infty(\theta_1,0,\ldots,0) = \bar{a}_0 + j\bar{a}_1\theta_1+ \bar{a}_2\theta_1^2 + \cdots .
\end{equation}  
Note that if $A$ fulfills Assumption~\ref{ass:relative}, then $\bar{a}_0 = 0$. If $A$ fulfills Assumption~\ref{ass:symmetry}, $\hat{a}_\infty(\theta)$ is real-valued and $\bar{a}_{1,3,\ldots} = 0$.

\subsection{Proof of Lemma~\ref{thm:riemannapprox}}
{Given that $\hat{\mathcal{A}}_\infty(\theta) $ is Hurwitz for $\theta \neq 0$, the \hnnorm density $\hat{p}(\theta)$ is continuous and bounded over the compact domain given by $\delta \le |\theta_i| \le \pi $ for $i = 1,\ldots,d$, and any fixed $\delta >0$. It is therefore Riemann integrable on that domain.} %In particular, $\hat{p}(\theta)$ is Riemann integrable on the domain $\delta \le |\theta_i| \le \pi $ for $i = 1,\ldots,d$ and any fixed $\delta >0$. }

{On the interval $\Delta < |\theta_i| < \delta $, allowing for $\Delta \rightarrow 0$, $\hat{p}(\theta)$  will instead be monotonic.  For simplicity, we show this through the scalar case in which $\hat{\mathcal{A}}_\infty(\theta) = \hat{a}_\infty(\theta)  = \sum_{k \in \mathbb{Z}^d} a_k e^{-j \theta \cdot k}$, which is negative for $\theta \neq 0$ as $\hat{\mathcal{A}}_\infty(\theta) $ is Hurwitz. Solving the Lyapunov equation \eqref{eq:lyap} then gives $\hat{P}(\theta) = \hat{p}(\theta) = -1/{2 \sum_{k \in \mathbb{Z}^d} a_k \cos(\theta \cdot k) }$ for $\theta \in \mathcal{R}^d \backslash \{0\}$. %Given the locality asssumption \eqref{eq:localitycondition}, its derivative in each coordinate direction % $i = 1,\ldots,d$ is
Its derivative in each coordinate direction $i = 1,\ldots,d$ is
\begin{equation}
\label{eq:derivativeh2density}
\frac{\mathrm{d}\hat{p}(\theta)}{\mathrm{d}\theta_i}  =  \frac{-2 \sum_{k \in \mathbb{Z}^d} a_k k_i\sin(k_1\theta_1 + \cdots k_d\theta_d)}{(2 \sum_{k \in \mathbb{Z}^d} a_k \cos(k_1\theta_1 + \cdots k_d\theta_d))^2}  .
\end{equation}
Now, note that $\mathrm{sgn}(\sin(kx)) = \mathrm{sgn}(x)$ for $\mathrm{for } |x| \le \frac{\pi}{k}$. Therefore, by the locality asssumption~\eqref{eq:localitycondition}, the derivative~\eqref{eq:derivativeh2density} satisfies $\frac{\mathrm{d}\hat{p}(\theta)}{\mathrm{d}\theta_i}  < 0$ for $\theta_i \in (0,\delta)$ and  $\frac{\mathrm{d}\hat{p}(\theta)}{\mathrm{d}\theta_i}  > 0$ for $\theta_i \in (-\delta,0)$ with $\delta \ge \pi/q$. 
%$\mathrm{sgn}(\sin(kx)) = \mathrm{sgn}(x) ~~\mathrm{for } |x| \le \frac{\pi}{k}. $  
The \hnnorm density $\hat{p}(\theta)$ is thus monotonically decreasing away from zero for $|\theta_i| \le \delta$, where $\delta$ can always be fixed. 
}
{A similar argument can be construed for when $\hat{\mathcal{A}}_\infty(\theta)$ is matrix-valued, in which case one considers matrix-valued coefficients of the $Z$-transform.}

{It is well-known that integrals of monotonic functions $f(x)$ can be estimated by upper and lower Riemann sums according to: $\int_{m}^{n+1}f(x)\mathrm{d}x \le \sum_{k = m}^n f(k) \le \int_{m-1}^{n}f(x)\mathrm{d}x$ if $f(x)$ decreasing (and vice versa if $f(x)$ increasing). We use this to bound the \emph{monotonic} part of the sum in \eqref{eq:samplingsum}: 
\begin{equation} \label{eq:Vn1}
\Vn^\delta = \frac{1}{L_\delta^d} ~\sum_{\mathclap{\substack{ \theta = \frac{2\pi}{L}n \\ |n_i|< \delta\frac{ L}{2 \pi} } } } \mathrm{tr} \left( \hat{\mathcal{B}}^*_\infty (\theta) \hat{P}(\theta) \hat{\mathcal{B}}_\infty(\theta) \right) \end{equation}
by the integral from $\Delta$ to $\delta$: $S^\delta(\Delta) := \int_{\Delta \le |\theta_1| \le \delta} \cdots  \int_{\Delta\le |\theta_d| \le \delta}   \hat{p}(\theta) \mathrm{d}\theta_1 \cdots \mathrm{d}\theta_d $
as 
$S^\delta\left(\frac{4\pi}{L}\right) \le  \Vn^\delta \le S^\delta\left(\frac{2\pi}{L}\right), $
since $2\pi/L$ and $4\pi/L$ are the first two wavenumbers, or sampling points in the sum. Here, $L_\delta$ is the number of summands for which $|n_i| < \delta\frac{ L}{2 \pi} $, corresponding to the domain where $\xtr$ is known to be monotonic. }

{For the remainder of the sum, we use the Riemann integrability away from zero. That is, let 
\[
\Vn^\pi = \frac{1}{(L-L_\delta)^d} \sum_{\mathclap{\substack{ \theta = \frac{2\pi}{L}n \\ |n_i|\ge \delta  \frac{L}{2 \pi} } } } \mathrm{tr} \left( \hat{\mathcal{B}}^*_\infty (\theta) \hat{P}(\theta) \hat{\mathcal{B}}_\infty(\theta) \right)
\] 
and note that $\lim_{L \rightarrow \infty} \Vn^\pi = S^\pi $, where $S^\pi:= \int_{\delta \le |\theta_1| \le \pi} \cdots  \int_{\delta\le |\theta_d| \le \pi}   \hat{p}(\theta)\mathrm{d}\theta_1 \cdots \mathrm{d}\theta_d $. That is, the sum converges to the integral. Therefore, at some $\bar{L}$, we will have that $|V_{\bar{N}}^\pi - S^\pi| < S^\delta(\frac{2\pi}{\bar{L}}) - S^\delta(\frac{4\pi}{\bar{L}})$, so that 
\[ S^\delta\left(\frac{4\pi}{L} \right) + S^\pi \le  \Vn^\delta + \Vn^\pi \le S^\delta\left(\frac{2\pi}{L}\right) + S^\pi, \]
for all $L \ge \bar{L}$, or $N > \bar{N}$, which is precisely equivalent to the statement of Lemma~\ref{thm:riemannapprox}.}

\subsection{Proof of Theorem~\ref{thm:newstabilitycons}}
Each matrix $\hat{\mathcal{A}}_\infty(\theta) = \begin{bmatrix}
\hat{a}_\infty(\theta) & \hat{b}_\infty(\theta) \\1 & \hat{f}_\infty(\theta)
\end{bmatrix} $ has eigenvalues
\(
\lambda_{1,2}= \frac{\hat{f}+ \hat{a}}{2} \pm \sqrt{ \left(\frac{\hat{f} - \hat{a}}{2}\right)^2 + \hat{b}  },
\)
where we omit the $\infty-$subscript and the argument $\theta$ of the individual Fourier symbols for notational compactness. 
The system is input-output stable if and only if $\mathrm{Re}\{ \lambda_{1,2}(\theta) \} <0$ for every $\theta \neq 0$ by Corollary~\ref{thm:stabilitycor}. To find necessary conditions for stability, it suffices to study this condition along one of the coordinate directions, so we let $\theta = (\theta_1,0,\ldots,0)$. 

A necessary condition for stability then becomes that
\begin{equation}
\label{eq:realpartcondition}
 \left| \mathrm{Re} \left\lbrace \frac{\hat{f}+ \hat{a}}{2}  \right\rbrace \right| > \left| \mathrm{Re}\left\lbrace \sqrt{ \left(\frac{\hat{f} - \hat{a}}{2}\right)^2 + \hat{b}  } \right\rbrace \right|,
\end{equation}
 for all $\theta_1 \in [-\pi,\pi]\backslash \{0\}$. If $\hat{a}, \hat{b}, \hat{f}$ are real-valued, \eqref{eq:realpartcondition} holds as long as $\hat{b} < \hat{a}\hat{f}$, which is true for example if $\hat{b}, \hat{a}, \hat{f} <0$. 

 If $\hat{a}, \hat{b}, \hat{f}$ are not all real-valued, the radicand on the right hand side (RHS) of \eqref{eq:realpartcondition} will be complex valued. Recall that for any complex number $z = |z|e^{j\phi}$, $\mathrm{arg}\{\sqrt{z}\} = \frac{1}{2} \mathrm{arg}\{z\} = \frac{1}{2}\phi$. This in particular means that if the argument $\phi$ is near $\pm \pi/2$, then $\mathrm{Re}\{\sqrt{z}\} =\sqrt{|z|}\cos(\phi/2) $ becomes large compared to $\mathrm{Re}\{z\} = |z| \cos \phi$. Here, this implies that \eqref{eq:realpartcondition} can only be satisfied if the imaginary part of the RHS radicand does not become ``too large'' compared to the real part.
 %this means that in order to satisfy \eqref{eq:realpartcondition}, the imaginary part of the radicand on the RHS cannot become ``too large'' compared to the real part.

We therefore study the radicand (now denoted $R$) on the RHS of \eqref{eq:realpartcondition} near $\theta_1=0$ by expanding it with the first terms of the Maclaurin expansions of the $Z$-transforms $\hat{a}, \hat{b}, \hat{f}$ as in~\eqref{eq:maclaurin}. Recalling that $B,F$ satisfy Assumption~\ref{ass:relative}, we get:

\vspace{-5mm}
\begin{small}
\begin{multline}
\label{eq:radicand}
R = :\left(\frac{\hat{f} - \hat{a}}{2}\right)^2 +\hat{b} \approx  \\  \frac{\bar{a}_0^2}{4} + \frac{(\bar{a}_2 - \bar{f}_2)^2}{4} \theta_1^4  + \left( \bar{b}_2- \frac{(\bar{a}_1-\bar{f}_1)^2}{4}  + \frac{\bar{a}_0(\bar{a}_2 - \bar{f}_2)}{2} \right) \theta_1^2 \\ +j\left[ \frac{(\bar{a}_1 - \bar{f}_1)(\bar{a}_2 - \bar{f}_2)}{2}\theta_1^3 + \left( \bar{b}_1 +\frac{\bar{a}_0(\bar{a}_1 - \bar{f}_1) }{2} \right) \theta_1 \right].
\end{multline}
\end{small}
Now, note that if $\bar{b}_1 \neq 0$, then %the imaginary part of $R$, 
$\mathrm{Im}\{R\}$, is \emph{linear} in~$\theta_1$ near $\theta_1= 0 $. Unless $\bar{a}_0 \neq 0$, $\mathrm{Re}\{R\}$ on the other hand, scales \emph{quadratically} in $\theta_1$, causing $\mathrm{Im}\{R\}$ to become arbitrarily many times larger than $\mathrm{Re}\{R\}$ as $\theta_1 \rightarrow 0$. This puts $\mathrm{arg}\{R\}$ near $\pm \pi/2$ and $\mathrm{Re}\{ \sqrt{R}\}$ becomes arbitrarily many times larger than $\mathrm{Re}\{R\}$. Since the LHS of \eqref{eq:realpartcondition} also grows quadratically in $\theta_1$ if $\bar{a}_0  = 0$, i.e., 
$  \left| \mathrm{Re} \left\lbrace \frac{\hat{f}+ \hat{a}}{2}  \right\rbrace \right|  \approx \left|  \frac{\bar{a}_0}{2} +\frac{\bar{a}_2 + \bar{f}_2}{2}\theta_1^2  \right|, $
it will also be smaller than $\mathrm{Re}\{ \sqrt{R}\}$ near $\theta_1= 0 $. We conclude that \eqref{eq:realpartcondition} cannot be fulfilled for all $\theta_1\in [-\pi,\pi]\backslash \{0\}$ in the case where $\bar{b}_1 \neq 0$, $\bar{a}_0 =0$. 

Necessary conditions for admissibility are therefore that $\bar{b}_1 = 0$, i.e., $\hat{b}$ real ($B$ symmetric), \emph{or} that $\bar{a}_0 \neq 0$ (absolute feedback in $A$).

\subsection{Proof of Theorem~\ref{thm:stabilityvehicle1}} 
The characteristic polynomial of the matrix $\hat{\mathcal{A}}_\infty(\theta) $
is 
\begin{equation}
\label{eq:chareqnvehicle}
p(\lambda, \theta) = \lambda^3 - (\hat{a} + \hat{g})\lambda^2 + (\hat{a}\hat{g} - \hat{f} - \hat{c})\lambda + \hat{a}\hat{f} - \hat{b},
\end{equation}
where we have again omitted the $\infty-$subscript and the argument $\theta$ of the Fourier symbols. %As in the previous proof, we study the system along $\theta = (\theta_1,0,\ldots,0)$.
Recall that all Fourier symbols are now real by Assumption~\ref{ass:symmetry}. We can therefore use the Routh-Hurwitz stability criteria which state that; given a characteristic polynomial $p(\lambda) = m_3\lambda^3 + m_2 \lambda^2 + m_1 \lambda +m_0,$ then necessary and sufficient criteria for stability are that \linebreak (i) $m_i >0$, $~i = 0,1,2,3$, and (ii) $m_2m_1 > m_3m_0$.

In the case of \eqref{eq:chareqnvehicle}, a necessary condition for satisfying (i) is that we do \emph{not} have $\hat{a} = 0$, $\hat{b} = 0$ simultaneously. {I.e., if $B = 0$, then we must have $A \neq 0$.} Otherwise, {the condition}~(i) can easily be satisfied, e.g. by ensuring $\hat{a}, \hat{b}, \hat{c}, \hat{f}, \hat{g}<0$. Assuming (i) is satisfied, consider (ii), which says that:
\(
 - (\hat{a} + \hat{g})(\hat{a}\hat{g} - \hat{f} - \hat{c})  > \hat{a}\hat{f} - \hat{b}.
\)
First, we note that if $\hat{b} = 0$, then this reduces to $-\hat{a}^2\hat{g} - \hat{a}\hat{g}^2 + \hat{a}\hat{c} +\hat{g}\hat{f} + \hat{g}\hat{c}>0$, which is also satisfied if $\hat{a}, \hat{c}, \hat{f}, \hat{g}<0$. 
For the case where $\hat{b} \neq 0$, we follow the approach in the previous proof and expand the inequality with the first terms of the Maclaurin expansions along $\theta_1$:
\vspace{-2mm}
\begin{multline}
\label{eq:routhabsolute}
- ( \bar{a}_0  +   (\bar{a}_2 + \bar{g}_2)\theta_1^2 ) ( \bar{a}_2\bar{g}_2\theta_1^4 + (\bar{a}_0\bar{g}_2 - \bar{f}_2 - \bar{c}_2 ) \theta_1^2  )\\
>  -\! \bar{b}_2\theta_1^2 + \bar{a}_0\bar{f}_2   \theta_1^2+\bar{a}_2\bar{f}_2 \theta_1^4
\end{multline} 
{Both sides of this inequality are positive if condition (i) above is satisfied.} Now, if the RHS of~\eqref{eq:routhabsolute} scales in lower powers of $\theta_1$ than the LHS, then near $\theta_1 = 0$ it becomes arbitrarily many times larger than the LHS, and \eqref{eq:routhabsolute} cannot be satisfied. In particular, if $\bar{b}_2  \neq 0$, then the RHS scales as $\theta_1^2$, 
{which is only true for the LHS if $\bar{a}_0 \neq 0$.} This concludes the proof.

\subsection{Proof of Lemma~\ref{lem:qscalingcons} } %changed April 1 because lemma is changed, changed April 6 because of mistake with imaginary part
To prove Lemma~\ref{lem:qscalingcons} we treat the two admissible feedback configurations given by Theorem~\ref{thm:newstabilitycons} separately. 
\subsubsection*{Case a) $B$ symmetric}
If $\hat{b}_\infty(\theta) $ is real, then
\begin{equation}
\label{eq:q1}
\varphi^c = \frac{\hat{b} \mathrm{Re}\{\hat{a}\}( \mathrm{Re}\{\hat{a}\} + \hat{f}) }{\hat{b}\hat{f} + \mathrm{Re}\{\hat{a}\}( \hat{b} -\mathrm{Im}\{\hat{a}\}^2 - (\mathrm{Re}\{\hat{a}\} + \hat{f})^2 )}.
\end{equation}
%where we have again omitted the $\infty$-subscript and the $\theta$-argument of the individual Fourier symbols. 
We notice immediately, that if $\hat{a}\equiv 0$, i.e., if $A = 0$, then $\varphi^c \equiv 0$, and $\hat{f} + \varphi^c$ scales just as $\hat{f}$.

Otherwise, recall that $\hat{f} \sim -\beta (\theta_1^2 + \ldots + \theta_d^2 )$ (for short: $\hat{f} \sim -\beta \theta^2$) by Lemma~\ref{prop:Fscale}. 
$B$ now has the same properties as $F$, so $\hat{b} \sim -\beta \theta^2$. $A$ on the other hand, may be asymmetric and have absolute feedback. Therefore, 
we in general have $\sum_{k \in \mathbb{Z}_L^d }a_k  = \hat{a}_0$, where $\hat{a}_0 \le 0$, and in line with~\eqref{eq:fouriercos} we obtain $\mathrm{Re}\{\hat{a}\} = \hat{a}_0 - \sum_{k \in \mathbb{Z}^d}a_k(1- \cos(\theta \cdot k ))$, so $\mathrm{Re}(\hat{a}) \sim \hat{a}_0 - \beta \theta^2$.
%in general, the condition \eqref{eq:relativedef} is not fulfilled and  $\sum_{k \in \mathbb{Z}_L^d }a_k : = \hat{a}_0$, where $\hat{a}_0 \le 0$ to guarantee stability. By \eqref{eq:fouriercos}, we then have that $\mathrm{Re}\{\hat{a}\} = \hat{a}_0 - \sum_{k \in \mathbb{Z}^d}a_k(1- \cos(\theta \cdot k )).$ Since the term $ \sum_{k \in \mathbb{Z}^d}a_k(1- \cos(\theta \cdot k ))$ has the same structure as $\hat{f}$, and we can write $\mathrm{Re}(\hat{a}) \sim \hat{a}_0 + \beta \theta^2$.
If $A$ is asymmetric, the imaginary part of its Fourier symbol is $\mathrm{Im}\{\hat{a}\} = - \sum_{k \in \mathbb{Z}^d}a_k \sin (\theta \cdot k)$. Through similar calculations as in the proof of Lemma~\ref{prop:Fscale}, we can derive the bound
$\mathrm{Im}\{\hat{a}\} ^2 \le (2(2q)^d +1)\sum_{k \in \mathbb{Z}^d}a_k^2 \sin^2(\theta \cdot k) = ((2q)^d +\frac{1}{2})\sum_{k \in \mathbb{Z}^d}a_k^2 \left(1 -  \cos(2 \theta \cdot k) \right)\le ((2q)^d +\frac{1}{2})(2q)^{(d+2)} ||a||_\infty^2 (\theta_1^2 + \cdots +\theta_d^2) . $
We can thus write $\mathrm{Im}\{\hat{a}\} ^2  \le \bar{c}_a\beta \theta^2$ with $\bar{c}\ge 0$. Clearly, it also holds $\mathrm{Im}\{\hat{a}\} ^2\ge 0$.

Now, consider the terms $\hat{b} - \left( \mathrm{Im}\{\hat{a}\} \right)^2$ in the denominator of~\eqref{eq:q1}. By the arguments in Appendix~\ref{sec:appScaling}, it holds $\hat{b} - \mathrm{Im}\{\hat{a}\} ^2 \sim -\beta\theta^2$. Inserting this, together with $\hat{f}, \hat{b}\sim -\beta \theta^2$, $\mathrm{Re}\{\hat{a}\} \sim \hat{a}_0 -\beta \theta^2$ into \eqref{eq:q1} gives
\[
\varphi^c \sim \beta \theta^2 \frac{- 2\beta\theta^2 + 2 \hat{a}_0 }{2\beta^2\theta^4 +  \beta(1 - 3\hat{a}_0)\theta^2 + \hat{a}_0^2}.
\]
This can be written as $\varphi^c \sim - \bar{\varphi} \beta\theta^2$, and we note that $\bar{\varphi}$ will be a bounded, positive constant for any $\beta$ and all  $\theta \in \mathcal{R}^d$. In fact, $\bar{\varphi} \le 2$ if  $\hat{a}_0 = 0$, or $\bar{\varphi}\le \frac{2}{|\hat{a}_0|}$ if $\hat{a}_0< 0$. 
Therefore, $\hat{f} + \varphi^c \sim -\beta\theta^2 - \bar{\varphi}\beta\theta^2 \sim -\beta \theta^2$, which is precisely \eqref{eq:qscaling}.

\subsubsection*{Case b) $B$ asymmetric}
If $B$ is not symmetric, we must by Theorem~\ref{thm:newstabilitycons} require $A$ to have %a sufficient amount of 
absolute feedback, so that $\hat{a}_\infty(\theta) \sim \hat{a}_0<0$. 
Inserting this into $\varphi^c$ gives
%Since the Fourier symbol of $A$ is then bounded by some scalar multiple of $\hat{a}_0$, it is not restrictive to assume $\hat{a}_\infty(\theta) = a_o$ for some $a_o<0$ from the beginning. We obtain:
\begin{equation}
\label{eq:q2}
%\varphi^c(\hat{a}, \hat{f}, \hat{b}) \sim \frac{ \hat{a}_0\mathrm{Re}\{\hat{b}\}(\hat{a}_0+\hat{f}) +\mathrm{Im}\{\hat{b}\}^2 }{\hat{a}_0(\mathrm{Re}\{\hat{b}\} - (\hat{a}_0+\hat{f})^2)+\hat{f}\mathrm{Re}\{\hat{b}\} }.
\varphi^c\sim \frac{ \hat{a}_0^2\mathrm{Re}\{\hat{b}\} +\mathrm{Im}\{\hat{b}\}^2 + \hat{a}_0\mathrm{Re}\{\hat{b}\}\hat{f}}{\hat{a}_0(\mathrm{Re}\{\hat{b}\} - (\hat{a}_0+\hat{f})^2)+\hat{f}\mathrm{Re}\{\hat{b}\} }.
\end{equation}
Now, $\mathrm{Im}\{\hat{b}\}^2 $ satisfies the {same inequality as $\mathrm{Im}\{\hat{a}\}^2$ above}. Since $\hat{f},\mathrm{Re}\{\hat{b}\} \sim -\beta \theta^2$, the numerator terms
%we have that the numerators terms 
$\hat{a}_0^2\mathrm{Re}\{ \hat{b} \} +\mathrm{Im}\{\hat{b}\}^2 \sim -\hat{a}_0^2\beta \theta^2$, provided that~$\hat{a}_0$ is sufficiently large to ensure admissibility. 
Inserting all scalings into \eqref{eq:q2} gives
\[\varphi^c \sim  \beta \theta^2 \frac{ \hat{a}_0^2 +\hat{a}_0}{\beta^2(\hat{a}_0 \!- \!1)\theta^4 + \beta \hat{a}_0(- 2\hat{a}_0 \! + \! 1)\theta^2 \! + \! \hat{a}_0^3} =:  -\bar{\varphi}\beta\theta^2\]
%\[\varphi^c(\hat{a}, \hat{f}, \hat{b}) \sim  \beta \theta^2 \frac{- a_o\beta\theta^2 -\beta + a_o^2 }{ \beta^2(a_o - 1)\theta^4 - \beta a_o( 2a_o - 1)\theta^2 + a_o^3} = -\tilde{c}\beta\theta^2.  \]
Here, $\bar{\varphi}$ is a positive constant, which for any $\beta$ and all $\theta \in \mathcal{R}^d$ satisfies $\bar{\varphi} \le \frac{1}{|\hat{a}_0|}$. We can again conclude that $\hat{f} + \varphi^c \sim -\beta\theta^2 $, which proves the lemma. 

\subsection{Proof of Lemma~\ref{thm:phifunction1} }
The function $\varphi^v$ in \eqref{eq:gramianrelative} is given as
\begin{equation}
\label{eq:phivfunction}
\varphi^v\! = \! \frac{ \hat{b}^2 \!+\! \hat{b} (\hat{a}\hat{c} \!+\! \hat{c}\hat{g} \! -\! \hat{a}\hat{f} \! -\! \hat{a}\hat{g}^2 \! -\! \hat{a}^2\hat{g})  \! -\! \hat{c}\hat{f}\hat{a}(\hat{a} \! +\! \hat{g})}{  \hat{b} \! -\! \hat{a}\hat{f} \!+\!\hat{a}^2(\hat{a} \!+\! \hat{g})  } 
\end{equation}
Now, the feedback operators $B,C,F,G$ have the relative measurement property~\eqref{eq:relativedef}, while $A$ need not to, so in line with Lemma~\ref{prop:fgvehicle}, we have $\hat{b}, \hat{c}, \hat{f}, \hat{g} \sim -\beta \theta^2$ and $\hat{a} \sim \hat{a}_0 - \beta\theta^2$ with $\hat{a}_0 \le 0$.
We consider the two cases given by the admissibility Theorem~\ref{thm:stabilityvehicle1} separately.
\subsubsection*{Case a) $B = 0$}
Substituting the scalings of the individual Fourier symbols into \eqref{eq:phivfunction} gives:
%\begin{small}
\[ 
\varphi^v\sim  \beta^2\theta^4 \frac{2 \beta \theta^2  - \hat{a}_0   }{ 2 \beta^2 \theta^4 +\beta (1 - 3 \hat{a}_0 ) \theta^2 + \hat{a}_0^2 } = \tilde{\varphi}\beta^2\theta^4.
\]
%\end{small}
For any $\beta$ and for all $\theta \in \mathcal{R}^d$, we identify $\tilde{\varphi}$ as a positive constant, with $\tilde{\varphi}\le \frac{1}{|\hat{a}_0|}$ if $\hat{a}_0 \neq 0$, $\tilde{\varphi}\le  2$ if $\hat{a}_0 = 0$. Therefore, $\hat{f}\hat{g} +\varphi^v \sim \beta^2\theta^4 + \tilde{\varphi}\beta^2\theta^4 \sim \beta^2\theta^4 $, which is precisely \eqref{eq:phiscalingvehicle}. 

\subsubsection*{Case b) $B \neq 0$}
If the operator $B$ is nonzero, $A$ is required by Theorem~\ref{thm:stabilityvehicle1} to have 
%a sufficient amount of 
absolute feedback, so $\hat{a}_0 < 0$. We can then set $\hat{a} \sim \hat{a}_0<0$ and:
\[\varphi^v \sim \beta^2 \theta^4 \frac{\beta (1- 2 \hat{a}_0) \theta^2 + 2 \hat{a}_0^2  -1 }{\beta \hat{a}_0 (\hat{a}_0+1) \theta^2 - \hat{a}_0^3} = \tilde{\varphi}\beta^2\theta^2.\]
Again, $\tilde{\varphi}$ can be identified as a bounded positive constant, so $\hat{f}\hat{g} + \varphi^v \sim \beta^2 \theta^4$ also in this case (provided $\hat{a}_0 \ge 1$, which signifies that the amount of absolute feedback in $A$ is sufficient to guarantee admissibility).

It remains to consider the case in which the feedback operator $C = 0$. This does not give a meaningful control design if $B = 0$, so it was not considered under case~a) above. Substituting $\hat{c} = 0$ and the scalings of remaining Fourier symbols into~\eqref{eq:phivfunction} gives
%\begin{small}
\[  \varphi^v \! \sim \!  \beta^2\theta^4 \!  \frac{ 2\beta^2\theta^4 \! - \! \beta( 3\hat{a}_0 \! + \! 1)\theta^2\!  + \! \hat{a}_0^2 \! + \! \hat{a}_0 \! -\!  1}{ 2\beta^3\theta^6 \! +\!  \beta^2 (1\! - \! 5\hat{a}_0)\theta^4\!  +\!  \beta(4\hat{a}_0^2 \!  -\!  \hat{a}_0 \! +\!  1)\theta^2 \! -\!  \hat{a}_0^3} ,\] %\end{small}
and the same conclusion as with $C \neq 0$ holds. %Notice that the case $C = 0$ is not a meaningful control design if $B = 0$, which is why it is not considered under case a). 

\subsection{Proof of Lemma~\ref{lem:controleffort}}
{Consider the dynamics~\eqref{eq:consensusdynamic}, but let the control signal $u = z + Fx$ be the output.}
We can {then} obtain the control signal variance in \eqref{eq:controlvariance} through the \hn norm from $w$ to $u$, divided by the total network size~$N$. We use the DFT~\eqref{eq:DFTdef} to block-diagonalize the system, and solve a Lyapunov equation for each wavenumber $n$. We obtain that 
\( \sum_{k \in \Znd} \mathbb{E}\{u_k^*u_k\} = \sum_{ n \in \Znd \backslash \{0\}  } \frac{\hat{b}_n - \hat{f}_n(\hat{f}_n + \hat{a}_n) }{2(\hat{a}_n + \hat{f}_n )} , \)
which is equivalent to 
\begin{equation}
\label{eq:l1norm}
N \mathbb{E} \{u_k^*u_k\} = \frac{1}{2}\left(||\hat{f}||_1 + ||\frac{\hat{b}}{\hat{a} + \hat{f}}||_1 \right) .
\end{equation}
The equivalence of the sum and the $l_1$-norm follows from the fact that we must have $\hat{f}_n, \hat{b}_n <0$ and $\hat{f}_n+\hat{a}_n<0$ for all $n$ in order to guarantee stability (see Theorem~\ref{thm:stability}).
Now, if $\hat{f}$ is the Fourier transform of a function array $f$, then $||\hat{f}||_\infty \le ||f||_1$ and $||f||_\infty \le \frac{1}{N}||\hat{f}||_1 $
%\begin{equation} \nonumber
%\label{eq:fourierbounds} ||\hat{f}||_\infty \le ||f||_1,~~||f||_\infty \le \frac{1}{N}||\hat{f}||_1
%\end{equation}
(see \cite{Bamieh2012}). Inserting in~\eqref{eq:l1norm} gives the first bound of the Lemma:
\(
||f||_\infty \le \frac{1}{N}||\hat{f}||_1 \le 2\mathbb{E} \{u_k^*u_k\}.
\)

It also holds that
\(
2N\mathbb{E} \{u_k^*u_k\} \ge ||\frac{\hat{b}}{\hat{a} + \hat{f}}||_1 \ge \frac{||\hat{b}||_1}{||\hat{a} + \hat{f}||_\infty  } \ge \frac{||\hat{b}||_1}{||\hat{a}||_\infty + ||\hat{f}||_\infty  } ,
\)
where the last equality follows from the triangle inequality. Now, we can use the fact that $||\hat{a}||_\infty \le ||\hat{a}||_1 \le (2q)^d||a||_\infty$ and substitute the bound {above }on $||\hat{f}||_\infty$ to get that
\(
2N\mathbb{E} \{u_k^*u_k\} \ge \frac{||\hat{b}||_1}{(2q)^d(||a||_\infty + 2\mathbb{E} \{u_k^*u_k\} ) } .
\)
Now, we use that $||b||_\infty \le \frac{1}{N}||\hat{b}||_1$ to rewrite this as
\[
4(2q)^d(\mathbb{E} \{u_k^*u_k\})^2 \!+ \! 2(2q)^d||a||_\infty \mathbb{E} \{u_k^*u_k\} \! -\!  ||b||_\infty \ge 0,
\]
which leads to the second bound of the Lemma. 

%%%%%%%%%%%%%%%%%%%%%%%%%%%%%%%%%%%%%%%%%%%%%%%%%%%%%%%%%%%%%%%%%%%%%%%%%%%%%%%%%%%
%%%%%%%%%%%%%%%%%%%%%%%%%%%%%%%%%%%%%%%%%%%%%%%%%%%%%%%%%%%%%%%%%%%%%%%%%%%%%%%%%%%

\section*{Acknowledgement}
We would like to thank Florian D\"orfler, Karl Henrik Johansson {and Bart Besselink} for a number of interesting discussions and insightful comments. {We are also grateful to the anonymous reviewers for their valuable feedback. }

\bibliographystyle{IEEETran}
\bibliography{emmasbib2015.bib,BassamBib.bib}

\begin{IEEEbiography}[{\includegraphics[width=1in,height=1.25in,clip,keepaspectratio]{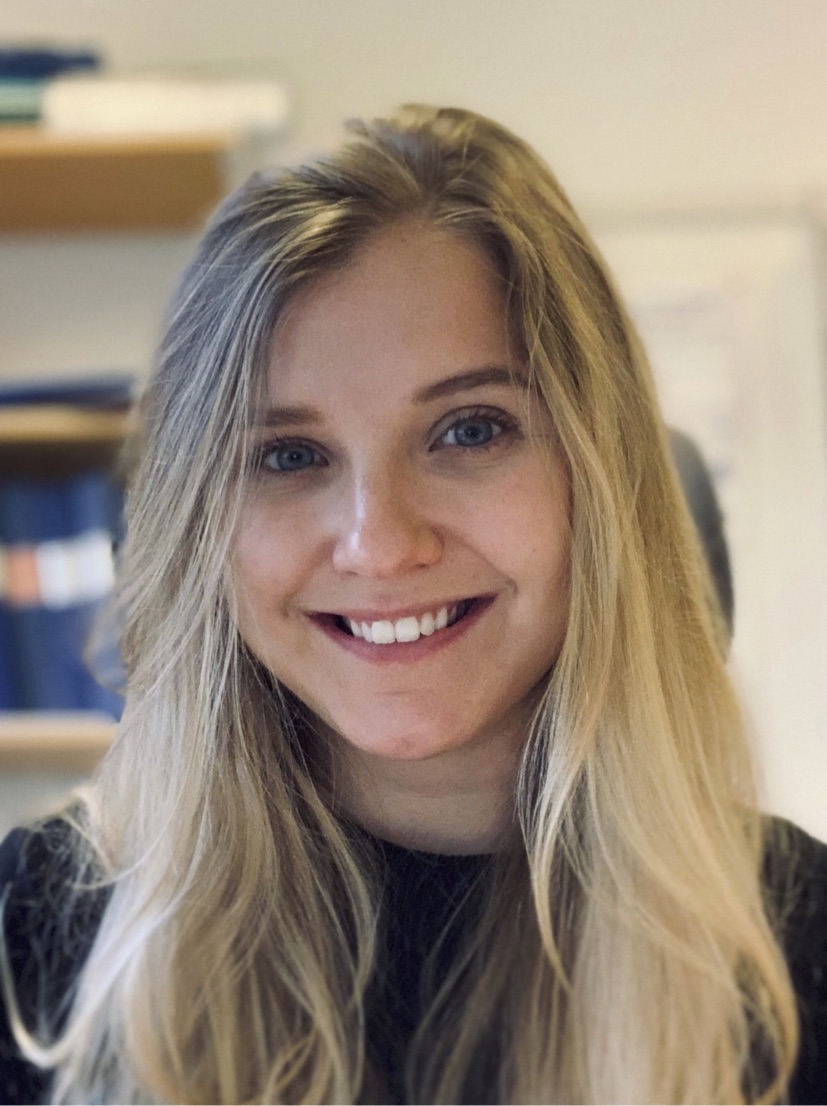}}]{Emma Tegling}
received her Ph.D. degree in Electrical Engineering in 2019 from KTH Royal Institute of Technology, Stockholm, Sweden, where she is currently a postdoctoral researcher at the Division of Decision and Control Systems. She received her M.Sc. degree in Engineering Physics in 2013, also from KTH. Dr. Tegling was a visiting researcher at Caltech in 2011, the Johns Hopkins University in 2013 and UC Santa Barbara in 2015. From 2013 to 2014 she was an analyst with Ericsson, Stockholm, Sweden. Her research interests are within analysis and control of large-scale networked systems.
\end{IEEEbiography}

\begin{IEEEbiography}[{\includegraphics[width=1in,height=1.25in,clip,keepaspectratio]{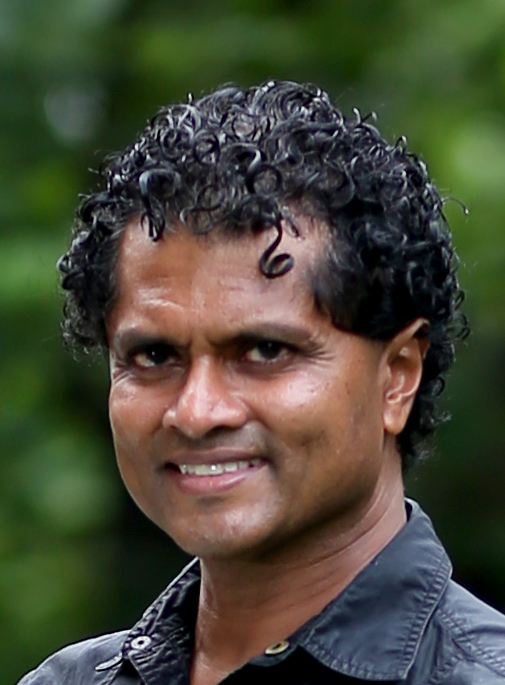}}]{Partha Mitra} received the Ph.D. degree in theoretical physics from Harvard University, Cambridge, MA, in 1993. He is Crick-Clay Professor of Biomathematics at Cold Spring Harbor Laboratory. He is also H N Mahabala Chair Professor (visiting) at IIT Madras and holds adjunct appointments at NYU Medical School and Weill Cornell Medical School. Dr Mitra is a member of the Theory Group at Bell Laboratories, Murray Hill (1993-2003), is a fellow of the American Physical Society, and a senior member of the IEEE.
\end{IEEEbiography}

\begin{IEEEbiography}[{\includegraphics[width=1in,height=1.25in,clip,keepaspectratio]{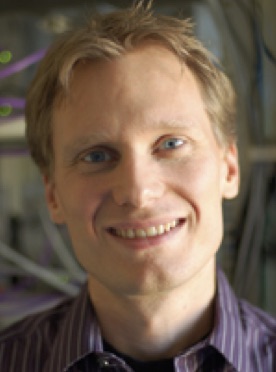}}]{Henrik Sandberg}
is Professor at the Division of Decision and Control Systems, KTH Royal Institute of Technology, Stockholm, Sweden. He received the M.Sc. degree in Engineering Physics and the Ph.D. degree in Automatic Control from Lund University, Lund, Sweden, in 1999 and 2004, respectively. His current research interests include security of cyberphysical systems, power systems, model reduction, and fundamental limitations in control. Dr. Sandberg received the Best Student Paper Award from the IEEE Conference on Decision and Control in 2004 and the Ingvar Carlsson Award from the Swedish Foundation for Strategic Research in 2007. He has served on the editorial board of IEEE Transactions on Automatic Control and is currently Associate Editor of the IFAC Journal Automatica.
\end{IEEEbiography}

\begin{IEEEbiography}[{\includegraphics[width=1in,height=1.25in,clip,keepaspectratio]{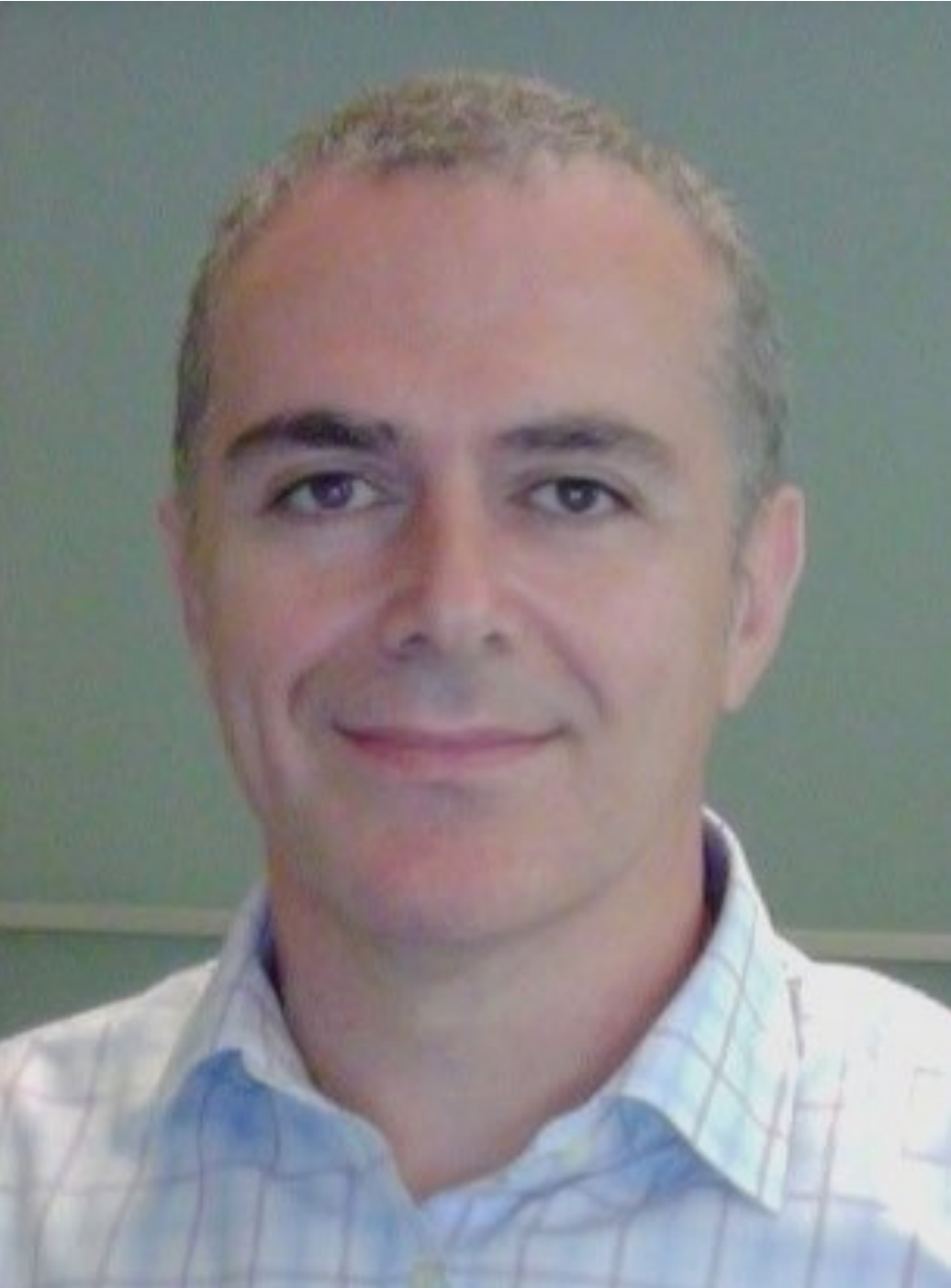}}]{Bassam Bamieh} is Professor of Mechanical Engineering at the University of California at Santa Barbara, CA. He received the B.S. degree in electrical engineering and physics from Valparaiso University, Valparaiso, IN, in 1983, and the M.Sc. and Ph.D. degrees from Rice University, Houston, TX, in 1986 and 1992 respectively.
%From 1991 to 1998 he was an Assistant Professor in the Department of Electrical and Computer Engineering and the Coordinated Science Laboratory of the University of Illinois at Urbana-Champaign. He is currently a Professor of Mechanical Engineering at the University of California at Santa Barbara. 
His current research interests are in optimal and robust control, distributed systems control, transition and turbulence modeling and control, and thermo-acoustic energy conversion devices.
Dr. Bamieh is a recipient of the AACC Hugo Schuck Best Paper Award, the IEEE CSS Axelby Outstanding Paper Award (twice), and an NSF CAREER Award. He is a Control Systems Society Distinguished Lecturer, a Fellow of IFAC, and a Fellow of the IEEE.
\end{IEEEbiography}

\end{document}